\documentclass[11pt, leqno]{article}
\usepackage{a4, indentfirst, latexsym,amssymb,amsmath}

  \newtheorem{Th}{Theorem}[section]
  \newtheorem{Prop}[Th]{Proposition}
  \newtheorem{Lem}[Th]{Lemma}
  \newtheorem{Cor}[Th]{Corollary}
  
  \newtheorem{Rem}[Th]{Remark}

  \def\pari#1#2{\makebox[10mm][r]{\em #1} \parbox[t]{138mm}{#2}}

  \def\Deg{\mathrm{Deg}}
  \def\la{\langle}
  \def\ra{\rangle}
  \def\cl{\overline }
  \def\Ker{\mathrm{Ker}}
  \def\Im{\mathrm{Im}}

  \def\d{\, d}

  \def\R{\mathbb{R}}

  \newcommand{\lma}{\lambda}
  \newcommand{\eps}{\varepsilon}

\begin{document}
\begin{center}
{\large{\bf Periodic solutions for nonlinear hyperbolic evolution systems}}

Aleksander \'{C}wiszewski\footnote{Corresponding author

}, Piotr Kokocki\\

{\em Faculty of Mathematics and Computer Science\\
Nicolaus Copernicus University\\
ul. Chopina 12/18, 87-100 Toru\'n, Poland }\\
\end{center}

\begin{abstract}
We shall deal with the periodic problem for
nonlinear perturbations of abstract hyperbolic evolution equations generating an evolution system of contractions. We prove an averaging principle for the translation along trajectories operator associated to the nonlinear evolution system, expressed in terms of the topological degree. The abstract results shall be applied to the damped hyperbolic partial differential equation.
\end{abstract}

\section{Introduction}

We shall be concerned with $T$-periodic solutions of the nonlinear evolution equation
$$
\dot u (t) = A(t)u(t) + F(t,u(t)), \quad t\in [0,T] \leqno{(P)}
$$
where $T>0$ is fixed, $\{A(t)\}_{t\in [0,T]}$ is a family of linear operators on a separable Banach space $E$ satisfying the so-called hyperbolic conditions and $F:[0,T]\times E\to E$ is
a $T$-periodic in time continuous map satisfying the local Lipschitz condition with respect to the second variable and having sublinear growth,
uniformly with respect to time. Moreover, it is also assumed that there is $\omega>0$ such that
$$
\|S_{A(t)}(s)\|\leq e^{-\omega s} \quad \mbox{ for \ } t\in [0,T] \mbox{ and } s\geq 0,
$$
where $S_{A(t)}$ stands for the $C_0$ semigroup generated by the operator $A(t)$,
and that there is $k\in [0,\omega)$ such that
$$
\beta (F([0,T]\times Q)) \leq k \beta (Q) \quad \mbox{ for any bounded } Q\subset E,
$$
where $\beta$ denotes the Hausdorff measure of noncompactness.
Under these assumptions, the translation along trajectories operators
$\Phi_t:E\to E$, $t\in [0,T]$, given by $\Phi_t(x):=u(t;x)$, $x\in E$, where $u(\cdot;x)$ stands for the solution of  $(P)$ with the initial condition $u(0)=x$, are well-defined and continuous. Moreover, for $t\in [0,T]$, one has $\beta(\Phi_t(Q))\leq e^{-(\omega-k)t} \beta (Q)$ for any bounded $Q\subset E$. This enables us to consider the topological degree of $I-\Phi_t$ and search $T$-periodic solutions corresponding to the fixed points of $\Phi_T$.
Our approach is based on the averaging idea, which says that if increasing the frequency in $(P)$, i.e. considering equations  $\dot u(t) = A(t/\lma)u(t) +F(t/\lma,u(t))$ with $\lma\to 0^+$,
then their solutions converge to solutions of the averaged equation
$$
\dot u(t) = \widehat A u(t) + \widehat F (u(t))
$$
where $\widehat A + \widehat F$ is the time averaged right-hand
side of $(P)$ (the precise meaning is explained in the sequel --
see Theorem \ref{18092008-1450}). Therefore, after rescaling time,
we study $T$-periodic solutions of equations
$$
\dot u(t) =  \lma A(t) u(t) + \lma F(t,u(t)), \quad t\in [0,T] \leqno{(P_\lma)}
$$
by means of the associated translation along trajectories operator $\Phi_{T}^{(\lma)}$.
We prove that, for small $\lma>0$, the topological degree of $I-\Phi_{T}^{(\lma)}$,
with respect to a proper open bounded $U\subset E$, is equal to the topological degree
$\Deg(\widehat A + \widehat F, U)$ -- see Theorem \ref{18092008-2220}.
This formula will imply the existence of $T$-periodic solutions provided $\Deg(\widehat A+\widehat F,U)\neq 0$.
In some natural cases the geometry of the right-hand side allows concluding the
nontriviality of the topological degree and get some a priori bounds estimates,
which provide effective criteria for the existence of $T$-periodic solutions  -- see Theorem \ref{24102008-0927}.\\
\indent The abstract hyperbolic type linear or semi-linear systems
and their applications to partial differential equations were
developed by Kato (see e.g. \cite{Kato}) and Tanabe (see e.g.
\cite{Tanabe} and references therein). Some existence results for
initial value problems associated with nonlinear perturbations of
evolution systems are standard and can be found e.g. in
\cite{Pazy}. As we need the continuity of translation along
trajectories and some related homotopies, we have to verify the
continuity and compactness of solutions as functions of initial
data and parameters. Moreover, due to some infinitesimal passages
related to the averaging method used in the paper, a parameterized
version of the representation formula must be derived. As a tool
we use the topological degree for so called $k$-set contraction
vector fields due to Sadovskii (see \cite{Akhmerov-etal} and
references therein) and Nussbaum (see \cite{Nussbaum}). The
topological degree for maps of the form $A+F$, where an invertible
operator $A$ generates a $C_0$ semigroup and $F$ is a continuous
$k$-set contraction is obtained as the degree of vector field
$I+A^{-1} F$, which is a standard -- see e.g. \cite{Kartsatos} and
some comments on the specific properties that we use are in
\cite{Cwiszewski-Kokocki}. Averaging methods combined with
topological degree and fixed point index were used in
\cite{Furi-Pera} to find periodic solutions for time dependent
vector fields on finite dimensional manifolds. Analogues of this
method, in the case of infinite dimensional Banach spaces was
stated in \cite{Cwiszewski-1}, \cite{Cwiszewski-2}, where periodic
solutions for the equations of the form $\dot u(t) = Au(t) +
F(t,u(t))$, with $A$ generating compact semigroups, were derived.
Also averaging methods together with Rybakowski's version of the
Conley index were used in \cite{Prizzi}, where the existence of
so-called recurrent solutions is studied for nonautonomous
parabolic equations. Periodic solutions for nonautonomous damped
hyperbolic equations has been also thoroughly studied in
\cite{Ortega1} and \cite{Ortega}. The present paper is a
continuation of \cite{Cwiszewski-Kokocki} where the periodic
problem is considered in the case where $A$ generates a
$C_0$-semigroup
of strict contractions and $F$ is a perturbation, i.e. the situation applicable to damped hyperbolic equations.\\
\indent The paper is organized as follows. In Section 2, we prove
a parameterized version of the representation theorem, which is a
useful framework for limit passages concerned with evolution
systems at all and also those considered in the next sections.
Section 3 is devoted to the properties of the translation along
trajectories operator such as the existence, continuity with
respect to the parameter and compactness. In Section 4 we deal
with the main result of the paper, that is the averaging method
for periodic solutions of $(P)$. Section 5 provides an example of
application to second order hyperbolic partial differential
equations.

\section{General Representation Theorem}

We start with a parameterized version of Theorem 3.5 from \cite[Ch. 3]{Pazy}.

\begin{Th}\label{29082008-1821}
Let $L:(0,+\infty)\times [0,1]\to {\cal L} (E,E)$, where $E$ is a Banach space, be a mapping such that
\begin{equation}\label{29082008-1144}
\|L(\lma,\mu)\|\leq 1\quad\mbox{ for \ } \ \lma> 0, \ \mu\in [0,1]
\end{equation}
and there is a dense subspace $V$ of $E$ such that
\begin{equation}\label{29082008-1145}
\lim_{\lma\to 0^+,\mu \to \mu_0} \lma^{-1}(L(\lambda,\mu)v - v) = A^{(\mu_0)}v \quad\mbox{ for \ } \ v\in V, \ \mu_0 \in [0,1],
\end{equation}
where, for each $\mu\in [0,1]$, $A^{(\mu)}:D(A^{(\mu)})\to E$ is a linear operator
such that $V\subset D(A^{(\mu)})$ and $(a_\mu I - A^{(\mu)})V$ is dense in $E$ for some $a_{\mu} >0$.\\
Then\\
\pari{(i)}{ for any $\mu\in [0,1]$, the operator $A^{(\mu)}$ is closable and its closure
$\overline{A^{(\mu)}}$
generates a $C_0$ semigroup of contractions $\{S_{\overline{A^{(\mu)}}}(t):E\to E\}_{t\geq 0}$;}\\[2mm]
\pari{(ii)}{for any sequence of a positive integers $(k_n)$ and a sequence
$(\lambda_n)$ in $(0,+\infty)$ such that $k_n\to\infty$, $k_n\lambda_n\to t$ as $n\to +\infty$, for some $t\geq 0$,
and any $(\mu_n)$ in $[0,1]$ with $\mu_n\to \mu_0$,
\begin{equation}\label{29082008-1157}
\lim_{n\to\infty}L(\lambda_n,\mu_n)^{k_n}x = S_{\overline{A^{(\mu_0)}}}\, (t)x \quad \mbox{ for each } x\in E;
\end{equation}}
\pari{(iii)}{for sequences $(k_n)$, $(\lma_n)$, $(\mu_n)$ and $t\geq 0$ as in {\em (ii)}
\begin{equation}\label{29082008-2233}
\lim_{n\to\infty}\lambda_n(I+L(\lambda_n, \mu_n)+L( \lambda_n,
\mu_n)^2+\ldots +L( \lambda_n, \mu_n)^{k_n-1})x\to \int_0^t
S_{\overline{A^{(\mu_0)}}}\,(\tau) x  \d \tau
\end{equation}
for each $x\in E$.}
\end{Th}
In the proof we shall use the following two Lemmata.
\begin{Lem}\label{29082008-1833}{\em(see \cite[Ch. 3, Theorem 4.5]{Pazy})}
If $(A_n)_{n\ge 1}$ is a sequence of operators generating $C_0$
semigroups $\{S_{A_n}(t)\}_{t\geq 0}$, $n\geq 1$, and $A:V\to E$ is a linear operator, where $V$ is a dense subspace of $E$, with the following properties\\[1mm]
\pari{(a)}{ there are $M\geq 1$ and $\omega\in\R$ such that
$\|S_{A_n}(t)\|\le Me^{\omega t}$ for any $n\geq 1$;} \\[1mm]
\pari{(b)}{ for every $v\in V$, $A_n v \to Av$ as $n\to\infty$;}\\[1mm]
\pari{(c)}{ there exists $\mu_0 >\omega$ such that
$(\mu_0 I -A)V$ is dense in $E$,} \\[1mm]
then the closure $\overline A$ of $A$ generates a $C_0$ semigroup $\{S_{\cl A}(t)\}_{t\ge 0}$ such that $$\|S_{\cl A}(t)\|\le Me^{\omega t} \quad\mbox{ for \ } t\geq 0 $$
and
$$\lim_{n\to\infty} S_{A_n}(t)x = S_{\cl A}(t)x \quad\mbox{ for \ }t\geq 0, \ x\in E.$$
The above convergence is uniform with respect to $t$ from bounded intervals.
\end{Lem}
\begin{Lem}\label{29082008-2156}{\em (see \cite[Ch. 3, Corollary 5.2]{Pazy})}
If $T\in {\cal L}(E,E)$  and $\|T\|\le 1$, then for any integer $n\geq 0$
and $x\in E$
$$
\|e^{(T - I)n}x - T^n x\|\leq \sqrt{n}\|x - Tx\|.
$$
\end{Lem}
\noindent{\bf Proof of Theorem \ref{29082008-1821}.}
(i) Define $A_{\lma}^{(\mu)}:E\to E$ by $A_{\lambda}^{(\mu)}:=\lambda^{-1}(L(\lambda,\mu)-I)$
and for any $\lma>0$, $\mu\in [0,1]$
and $t\geq 0$, put $S_{\lambda}^{(\mu)} (t):=\exp(tA_{\lma}^{(\mu)})$. Clearly, in view of (\ref{29082008-1144}), for any $\lma>0$ and $\mu\in [0,1]$
\begin{equation}\label{16072009-1816}
\|S_{\lma}^{(\mu)}(t)\|\le e^{-t/\lma} \sum_{k=0}^\infty(t/\lma)^k\frac{\|L(\lma,\mu)^k\|}{k!}\leq
e^{-t/\lma}\sum_{k=0}^\infty\frac{(t/\lma)^k}{k!} = 1.
\end{equation}
If $\lma_n\to 0^+$ and $\mu_n\to \mu_0$, then due to (\ref{29082008-1145}),
$$
\lim_{n\to \infty} A_{\lma_n}^{(\mu_n)} v = A^{(\mu_0)}v \quad\mbox{ for \ }
v\in V.
$$
By the assumption, there is $a_{\mu_0}>0$ such that $(a_{\mu_0} I - A^{(\mu_0)})V$ is dense in $E$ and, in view of  Lemma \ref{29082008-1833},
we infer that $A^{(\mu_0)}$ is closable and
its closure $\overline{A^{(\mu_0)}}$ generates $C_0$ a semigroup $\left\{ S_{\overline{A^{(\mu_0)}}}\, (t)\right\}_{t\geq 0}$ of bounded linear operators on $E$, such that $\|S_{\overline{A^{(\mu_0)}}}\, (t)\|\leq 1$ for any $t\geq 0$ and furthermore
\begin{equation}\label{29082008-2222}
S_{\lma_n}^{(\mu_n)} (t) x\to S_{\overline{A^{(\mu_0)}}}\, (t)x \quad\mbox{ for any $x\in E$}, \mbox{ as } n\to \infty
\end{equation}
uniformly for $t$ from bounded subintervals of $[0,+\infty)$.\\
\indent (ii) Let the sequence of a positive integers $(k_n)$, the
sequence $(\lambda_n)$ in $(0,+\infty)$ and $(\mu_n)$ in $[0,1]$
be such that $k_n\to\infty$, $k_n\lambda_n\to t$ as $n\to
+\infty$, for some $t\geq 0$, and $\mu_n\to \mu_0$ as $n\to
+\infty$. Then for any $v\in V$ and $n\geq 1$
\begin{align}\label{29092008-2220}
& \| L(\lambda_n,\mu_n)^{k_n}v - S_{\overline{A^{(\mu_0)}}}\, (t)v \| \leq
\| L(\lambda_{n}, \mu_{n})^{k_n}v - S_{\lma_n}^{(\mu_n)}(\lma_n k_n)v\| \\ \nonumber
& \quad + \|S_{\lma_n}^{(\mu_n)}(\lma_n k_n)v-S_{\overline{A^{(\mu_0)}}}\, (\lma_n k_n)v\|
+\|S_{\overline{A^{(\mu_0)}}}\, (\lma_n k_n)v-S_{\overline{A^{(\mu_0)}}}\, (t) v\|\nonumber.
\end{align}
By Lemma  \ref{29082008-2156} and (\ref{29082008-1145}), for any $v\in V$
\begin{align}\nonumber
\|L(\lambda_n,\mu_n)^{k_n}v - S_{\lma_n}^{(\mu_n)} (\lma_n k_n)v \| & =
\| e^{k_n(L(\lambda_n,\mu_n)-I)}v-L(\lambda_n,\mu_n)^{k_n}v\| \\ \nonumber
& \le \sqrt{k_n}\|v - L(\lambda_n,\mu_n)v\| \\
& = \sqrt{\lma_n} \sqrt{k_n\lma_n} \|\lma_n^{-1}(v-L(\lambda_n, \mu_n)v)\|
\to 0 \quad\mbox{ as } n\to\infty.
\end{align}
Consequently, by (\ref{29082008-1144}), (\ref{16072009-1816}) and the density of $V$ in $E$, we obtain that
\begin{equation}\label{16082008-1832}
\|L(\lambda_n,\mu_n)^{k_n}x - S_{\lma_n}^{(\mu_n)} (\lma_n k_n)x \| \to 0 \quad \mbox{ for each } x\in E, \mbox{ as } n\to \infty .
\end{equation}
Furthermore, in view of the uniform convergence on bounded intervals in (\ref{29082008-2222}), one has
$$
\|S_{\lma_n}^{(\mu_n)}(\mu_n k_n)x-S_{\overline{A^{(\mu_0)}}}\, (\lma_n k_n)x \| \to 0 \quad \mbox{ for any } x\in E, \mbox{ as } n\to \infty.
$$
This, together with (\ref{29092008-2220}), (\ref{16082008-1832}) and the continuity of the semigroup $S_{\overline{A^{(\mu_0)}}}$, gives (\ref{29082008-1157}). \\
\indent (iii) Take any $v\in V$ and observe that
$$
\left\|\lma_n \sum_{k=0}^{k_n-1} L(\lambda_n,\mu_n)^k v  -
\int_{0}^{t} S_{\overline{A^{(\mu_0)}}}\,(\tau)v \d \tau
\right\|\leq I_{n}^{(1)} + I_{n}^{(2)} + I_{n}^{(3)},
$$
where
\begin{align*}
I_{n}^{(1)} & :=\left\|\lma_n \sum_{k=0}^{k_n-1} L(\lambda_n,\mu_n)^k v  -
\lma_n \sum_{k=0}^{k_n-1} S_{\lma_n}^{(\mu_n)}(k\lma_n)v\right\|,\\
I_{n}^{(2)} & :=\left\| \lma_n \sum_{k=0}^{k_n-1} S_{\lma_n}^{(\mu_n)}(k\lma_n)v -
\lma_n \sum_{k=0}^{k_n-1} S_{\overline{A^{(\mu_0)}}}\,(k\lma_n)v \right\|,\\
I_{n}^{(3)} & :=\left\| \lma_n \sum_{k=0}^{k_n-1}
S_{\overline{A^{(\mu_0)}}}\,(k\lma_n)v - \int_{0}^{t}
S_{\overline{A^{(\mu_0)}}}\, (\tau)v \d \tau\right\|.
\end{align*}
First, in view of Lemma \ref{29082008-2156} and (\ref{29082008-1145}), one has
\begin{align}\label{29082008-2320}
I_{n}^{(1)} & \leq k_n\lma_n \max\{
\|L(\lambda_n,\mu_n)^k-e^{k(L(\lma_n, \mu_n) - I)}\| \mid k=1,\ldots, k_n-1\}\nonumber\\
& \leq  k_n\lma_n \max\{\sqrt{k} \|v-L(\lma_n, \mu_n)v \| \mid k=1,\ldots, k_n-1\}\\
& \leq \sqrt{\lma_n} (k_n\lma_n)^{3/2} \|\lma_n^{-1}(v-L(\lma_n,\mu_n)v) \| \to 0 \quad\mbox{ as } n\to \infty\nonumber.
\end{align}
Furthermore, by the uniform convergence in (\ref{29082008-2222}) on the interval $[0,\overline t]$,
where $\overline t := \sup_{n\geq 1} k_n\lma_n$, we get
\begin{equation}\label{29082008-2321}
I_{n}^{(2)} \!\leq\! k_n\lma_n \max\{ \|S_{\lma_n}^{(\mu_n)}(\!k\lma_n\!)v \!-\! S_{\overline{A^{(\mu_0)}}}\,(\!k\lma_n\!)v\|\! \mid \!k=0,\ldots,\! k_n\!-\!1\}\!\to 0 \quad \mbox{ as } n\to\infty.
\end{equation}
It is also clear that $I_{n}^{(3)}\to 0$ as $n\to \infty$, which along with (\ref{29082008-2320}) and (\ref{29082008-2321}) implies that (\ref{29082008-2233}) is satisfied for $x\in V$.
Finally, since for any $n\geq 1$
$$
\left\|\lma_n\sum_{k=0}^{k_n-1} L(\lma_n,\mu_n)^k- \int_{0}^{t}
S_{\overline A^{(\mu_0)}}\,(\tau) \d\tau \right\|\leq k_n\lma_n +
t < C $$ for some constant $C>0$ independent of $n$ and since $V$
is dense in $E$, one has the required convergence for each $x\in
E$.  \hfill $\square$

\section{Continuity and compactness properties for solution \\ operator}

A family $\{R(t,s)\}_{0\leq s\leq t\leq T}$, $T>0$ of bounded linear operators on a Banach space $E$ is called an {\em evolution system} provided $R(t,t)=I$ for each $t\in [0,T]$,
$R(t,s)=R(t,r)R(r,s)$ if only $0\leq s\leq t\leq T$ and for any $x\in E$, the map
$(t,s)\mapsto R(t,s)x$ is continuous. A family $\{R^{(\lma)}\}_{\lma\in [0,1]}$ of evolution systems is called {\em continuous} if, for any $x\in E$ and $(\lma_n)$  in $[0,1]$ with $\lma_n\to \lma$, $R^{(\lma_n)}(t,s)x\to R^{(\lma)}(t,s)x$ uniformly with respect to $t,s\in[0,T]$ with $s\leq t$.\\
\indent Evolution systems are naturally determined by time-dependent families of linear operators. Namely,
if $\{ A (t) \}_{t\in [0,T]}$ is a family of linear operators on a Banach space $E$ such that for any $s\in [0,T]$ and $x\in E$, the problem
$$\left\{
\begin{array}{l}
\dot u (t) = A(t)u(t), \quad t\in [s,T]\\
u(s)=x
\end{array}\right.
$$
admits (in some sense) a unique solution $u_{s,x}:[s,T]\to E$,  then the corresponding evolution system $\{R(t,s)\}_{0\leq s\leq t\leq T}$ is given by $R(t,s)x:=u_{s,x} (t)$, for $t\in [s,T]$. A particular type of evolution systems -- the so-called {\em hyperbolic evolution systems}, will be discussed in details at the end of this section.
\begin{Prop}\label{30122008-1629}
Suppose that $\left\{R^{(\lma)}\right\}_{\lma\in [0,1]}$ is a continuous family of evolution systems and the operator $\Sigma:E\times L^1([0,T], E)\times [0,1] \to C([0,T],E)$ is given by
$$
\Sigma (x,w,\lma) (t):= R^{(\lma)}(t,0)x+ \int_{0}^{t} R^{(\lma)}(t,s) w(s) \d s.
$$
Then\\
\makebox[10mm][r]{\em (i)} \parbox[t]{138mm}{$\Sigma$ is continuous;}\\
\pari{(ii)}{if $K\subset E$ is relatively compact and $W\subset L^1([0,T],E)$ is such that there is $c\in L^1([0,T])$ with  $\|w(t)\|\leq c(t)$ for any $w\in W$ and a.e. $t\in [0,T]$, then
$\Sigma (K\times W\times [0,1])$ is relatively compact if and only if the set
$\{ u (t) \mid u\in \Sigma (K\times W\times [0,1])\}$ is relatively compact for any $t\in [0,T]$.}
\end{Prop}
\begin{Rem}\label{07012009-1059} {\em
\indent (a) Under the above notation, if $t,t+h\in [0,T]$ with $h>0$, then
$$
\Sigma (x,w,\lma) (t+h) = R^{(\lma)}(t+h,t)\Sigma(x,w,\lma)(t) +
\int_{t}^{t+h} R^{(\lma)} (t+h,s)w(s)\ \d s,
$$
which follows directly from the definition of $\Sigma$ and the properties of evolution systems.\\
\indent (b) If $\left\{R^{(\lma)}\right\}_{\lma\in [0,1]}$ is a continuous family of evolution systems, then for any $x\in E$ the set $\{ R^{(\lma)} (t,s)x \ | \ 0\leq s\leq t\leq T, \lma\in [0,1]\}$ is bounded. Hence, in view of the uniform boundedness principle, there exists $M\geq 0$ such that
$$
\|R^{(\lma)}(t,s)\| \leq M \quad\mbox{ for any } t,s\in [0,T] \mbox{ with } s\leq t\mbox{ and } \lma\in [0,1].
$$}
\end{Rem}
{\bf Proof of Proposition \ref{30122008-1629}.} (i) Let $x_n\to x_0$ in $E$, $w_n\to w_0$ in $L^1([0,T],E)$ and $\lma_n\to \lma_0$.
Clearly, by Remark \ref{07012009-1059} one has
\begin{align}\label{07012009-1130} \nonumber
\|R^{(\lambda_n)}(t,0)x_n\! - \! R^{(\lambda_0)}(t,0)x_0 \|
& \leq \| R^{(\lambda_n)}(t,0)x_n\!-\!R^{(\lambda_n)}(t,0)x_0 \| + \\
& \hspace{26mm} + \| R^{(\lambda_n)}(t,0)x_0-R^{(\lambda_0)}(t,0)x_0\| \\
& \leq M\|x_n-x_0\| + \|R^{(\lambda_n)}(t,0)x_0 - R^{(\lambda_0)}(t,0)x_0 \|\nonumber
\end{align}
and hence, by the continuity of the family $\{R^{(\lma)}\}_{\lma\in [0,1]}$, we infer that $\|R^{(\lambda_n)}(t,0)x_n - R^{(\lambda_0)}(t,0)x_0
\| \to 0$ as $n\to +\infty$, uniformly with respect to
$t\in [0,T]$. In a similar manner
\begin{align} \nonumber
\left\| \int_{0}^{t}\!\! R^{(\lma_n)}(t,s) w_n(s) \d s \!-\!\!\int_{0}^{t}\!\! R^{(\lma_0)}(t,s) w_0(s) \d s\right\| \\ \nonumber
& \hspace{-30mm}\leq \int_{0}^{t}\!\! \|R^{(\lma_n)}(t,s) w_n(s) \!-\! R^{(\lma_0)}(t,s) w_0(s)\|\d s \\
& \hspace{-30mm}\leq  \int_{0}^{t} \!\!\|R^{(\lma_n)}(t,s) w_n(s) \!-\! R^{(\lma_n)}(t,s) w_0(s)\| \d s \\ \nonumber
& \hspace{-30mm} \quad + \int_{0}^{t}\!\! \|R^{(\lma_n)}(t,s) w_0(s)\!-\! R^{(\lma_0)}(t,s) w_0(s)\| \d s\\ \nonumber
& \hspace{-30mm} \leq M \|w_n-w_0\|_{L^1([0,T],E)} + \int_{0}^{T} \varphi_n (s) \d s, \nonumber
\end{align}
where functions $\varphi_n:[0,T]\to \R$, $n\geq 1$, are given by
$$
\varphi_n (s) := \sup_{\sigma,\tau\in [0,T],\tau\geq \sigma} \| [ R^{(\lma_n)}(\tau,\sigma)- R^{(\lma_0)} (\tau,\sigma)] w_0(s) \|.
$$
It is easy to check that functions $\varphi_n$, $n\ge 1$ are
measurable and, by the continuity of $\{R^{(\lma)}\}_{\lma\in
[0,1]}$ and Remark \ref{07012009-1059}, we infer that, for a.e.
$s\in [0,T]$, $\varphi_n(s)\to 0$ as $n\to+\infty$. On the other
hand $0\leq \varphi_n (s) \leq 2M\|w_0 (s)\|$, for $s\in [0,T]$,
which in view of the Lebesgue dominated convergence theorem
gives $\int_{0}^{T} \varphi_n(s) \d s \to 0$ as $n\to + \infty$ and together with (\ref{07012009-1130}) proves (i).\\
\indent (ii) Suppose that the set $\{ u (t) \mid u\in \Sigma (K\times W\times [0,1])\}$ is relatively compact for any $t\in [0,T]$.
Take any $\eps>0$ and fix $t\in [0,T]$. Let $\delta > 0$ be such that
$\int\limits_{[t-\delta,t+\delta]\cap [0,T]} c(s) \d s <\eps/3 M$.
Suppose that $t\in [0,T)$. Since the set
$Q_t:=\overline{\{ \Sigma(x,w,\lma)(t)\mid x\in K, w\in W, \lma\in [0,1]\}}$
is compact, one may eventually decrease $\delta > 0$ so that
$$
\|R^{(\lma)}(t+h,t)z-z\|<\eps/2, \quad \mbox{ if $t+h \leq T$, $h\in[0,\delta)$, $\lambda\in[0,1]$, $z\in Q_t$.}
$$
Now take $h\in [0,\delta)$ such that $t+h\in [0,T]$. Then, denoting
$\Sigma:=\Sigma(x,w,\lma)$ for any $(x,w,\lma)\in K\times W\times [0,1]$, one has
\begin{align*}
\|\Sigma (t+h)-\Sigma(t)\| & \leq
\| \Sigma (t+h)- R^{(\lma)}(t+h,t)\Sigma (t)\|+\|R^{(\lma)}(t+h,t)\Sigma (t)-\Sigma(t)\|\\
& \leq \int_{t}^{t+h} \|R^{(\lma)}(t+h,s)w(s)\| \d s + \|R^{(\lma)}(t+h,t)\Sigma (t)-\Sigma(t)\|\\
& \leq M \int_{t}^{t+h} c(s) \d s + \eps/2 <\eps.
\end{align*}
If $t\in (0,T]$, then take any $\delta_1 \in (0, \min\{t,\delta\}]$.
Since the set
$$
Q_{t-\delta_1}:=\overline{\{ \Sigma(x,w,\lma)(t-\delta_1)\mid x\in K, w\in W, \lma\in [0,1]\}}
$$
is compact, there exists $\delta' \in (0, \delta_1]$ such that
$$\|R^{(\lambda)}(t-h,t-\delta_1)z - R^{(\lma)}(t,t-\delta_1)z\|\leq \eps/3 \quad \mbox{ for any $h\in [0,\delta')$, $\lambda\in [0,1]$, $z\in Q_{t-\delta_1}$}.$$
In consequence, for any $x\in K$, $w\in W$, $\lma\in [0,1]$ and $h\in [0,\delta')$
\begin{align*}
\|\Sigma (t-h)-\Sigma(t)\| & \leq \|\Sigma(t-h)-R^{(\lma)}(t-h,t-\delta_1)\Sigma(t-\delta_1) \| \\
& \quad + \| R^{(\lambda)}(t-h,t-\delta_1)\Sigma(t-\delta_1) - R^{(\lma)}(t,t-\delta_1)\Sigma(t-\delta_1)\| \\
& \qquad + \|R^{(\lma)}(t,t-\delta_1)\Sigma(t-\delta_1) - \Sigma(t)\|\\
& \leq  \left\|\int_{t-\delta_1}^{t-h} R^{(\lma)}(t-h,s) w(s) \d s \right\|+\eps/3+
\left\|\int_{t-\delta_1}^{t} R^{(\lma)}(t,s) w(s) \d s \right\|\\
& \leq M \int_{t-\delta_1}^{t-h} c(s) \d s + \eps/3 + M\int_{t-\delta_1}^{t} c(s) \d s
\leq \eps.
\end{align*}
Hence, the set $\{\Sigma(x,w,\lma) \}_{x\in K, w\in W, \lma\in [0,1]}$ is equicontinuous at any $t\in [0,T]$, which due to the Ascoli-Arzel\`{a} theorem, completes the proof of (ii). \hfill $\square$

We state basic continuity and compactness results for the translation along trajectories operator
for a perturbed (possibly nonlinear) equation
\begin{equation}\label{14092008-0948}
\left\{\begin{array}{l}
\dot u(t)=A(t)u(t)+F(t,u(t)), \quad t\in [t_0,T] \\
u(t_0) = x_0
\end{array}\right.
\end{equation}
where $\{ A(t) \}_{t\in [0,T]}$ is a family of linear operators on Banach space $E$ with the associated evolution system $\{R(t,s)\}_{0\leq s\leq t\leq T}$, $F:[0,T]\times E\to E$ is a continuous map, $t_0\in[0,T)$ and $x_0\in E$. Recall (after \cite{Pazy}) that a continuous function $u:[t_0,T]\to E$ is called a {\em mild solution } of (\ref{14092008-0948}) if and only if
$$
u(t)=R(t,t_0) x_0 + \int_{t_0}^{t} R(t,s)F(s,u(s)) \d s
\quad\mbox{ for any } t\in [t_0,T].
$$
For the purpose of next sections we consider below a parameterized framework.
\begin{Prop} \label{09012009-2320}
For each $\lma\in [0,1]$, let $\{A^{(\lma)}(t)\}_{t\in [0,T]}$ be
family of operators on a separable Banach space $E$ having
associated evolution system $R^{(\lma)}$. If the family
$\{R^{(\lma)}\}_{\lambda\in[0,1]}$ is continuous and
$F:[0,T]\times E\times [0,1]\to E$ is a continuous map being \\[1mm]
\noindent $(F_{1})_{par}$ \parbox[t]{138mm}{locally Lipschitz in the second variable uniformly with respect
to the other variables, i.e. for each $x\in E$, there exists $r_x>0$ and $L_x>0$ such that for each
$t\in [0,T]$, $x_1, \ x_2\in B(x,r_x)$ and $\lma \in [0,1]$
$$\|F(t,x_1, \lma) - F(t,x_2,\lma)\|\le L_x \|x_1 - x_2\|;$$}
\noindent $(F_{2})_{par}$ \parbox[t]{138mm}{ of  sublinear growth in the second variable uniformly with respect to the others, i.e. there is $c>0$ such that
$$
\|F(t,x,\lma)\|\le c(1 + \|x\|) \quad \mbox{ for any }  x\in E, t\in [0,T], \lma\in [0,1];
$$ } \\[1mm]
\noindent $(F_{3})_{par}$ \parbox[t]{138mm}{a $k$-set contraction, i.e. there exists $k\geq 0$ such that
 $$
 \beta(F([0,T]\times Q\times [0,1]))\le k\beta(Q)
\quad  \mbox{ for any bounded } Q\subset E,$$}\\[1em]
\indent then\\
\pari{(i)}{ {\em (Existence)} for any $x\in E$ and $\lambda \in [0,1]$,
the initial value problem
\begin{equation}\label{07012009-1523}
\left\{
\begin{array}{l}
\dot u(t) = A^{(\lma)}(t)u(t)+F(t,u(t),\lma),\quad t\in [0,T]\\
u(0)=x
\end{array}
\right.
\end{equation}
admits a unique mild solution $u(\,\cdot\, ;0,T,x,\lambda)$;}\\[2mm]
\pari{(ii)}{{\em (Continuity)} if $( x_n, \lambda_n )\to ( x_0 , \lambda_0 )$
in $E\times [0,1]$, then $$
u(\, \cdot \, ;0,T,x_n,\lambda_n) \to u(\,\cdot\, ;0,T,x_0,\lambda_0) \ \ \ \mbox{ in } C([0,T],E);
$$}
\pari{(iii)}{ {\em (Compactness)} if additionally, there is  $\omega >0$ such that
\begin{equation}\label{03082009-1958}
\|R^{(\lma)}(t,s) \| \leq e^{-\omega(t-s)} \quad \mbox{ for \ } 0
\leq s \leq t\leq T
\end{equation}
and, for any $t\in [0,T]$, $\Phi_t:E\times [0,1] \to E$ is given by $\Phi_t (x,\lambda):= u(t;0,T,x,\lambda),$
then, for any bounded $Q\subset E$ and $t\in [0,T]$, the set $\Phi_t (Q\times [0,1])$ is bounded and
\begin{equation*}
\beta \left( \Phi_t (Q\times [0,1]) \right) \leq  e^{(k-\omega) t} \beta(Q).
\end{equation*}}
\end{Prop}
\begin{Lem} \label{19092008-0029} {\em (see \cite{Deimling}, \cite{ObuKaZ})}
Suppose that $E$ is a separable Banach space, $W\subset L^1([a,b],E)$ is countable and there is $c\in L^1([a,b])$
such that $\|w(t)\|\leq c(t)$,  for all $w\in W$ and a.e. $t\in [a,b]$,
and let $\phi:[a,b]\to\R$ be given by $\phi (t):= \beta(\{w(t)\ | \ w\in W \})$.
Then $\phi\in L^1([a,b])$ and
$$
\beta \left(\left\{
\int_{a}^{b} w(\tau) \d \tau\,|\, w \in W\right\} \right)
\leq \int_{a}^{b} \phi(\tau) \d \tau.
$$
\end{Lem}
\begin{Lem}\label{08012009-2310}
Let $\{R^\lma\}_{\lma\in [0,1]}$ be a continuous family of evolution system such that
\begin{equation}\label{03082009-1907}
\|R^{(\lma)}(t,s) \| \leq M e^{\omega(t-s)} \quad \mbox{ for \ } 0
\leq s \leq t\leq T \mbox{ and }\lma\in [0,1]
\end{equation}
where $M>0$ and $\omega\in\R$ are constants. Then, for any bounded $Q\subset E$ and $s,t \in [0,T]$ with $s\leq t$,
$$
\beta \left( \{ R^{(\lma)}(t,s)x\mid x\in Q, \, \lma\in [0,1] \}\right) \leq M e^{\omega(t-s)} \beta(Q).
$$
\end{Lem}
The proof is analogical to the proof of Lemma 2.1 from \cite{Cwiszewski-Kokocki}.

\noindent {\bf Proof of Proposition \ref{09012009-2320}.}
(i) follows by standard arguments (see e.g. \cite{Pazy} and \cite{Daners}).\\
\indent The proofs of (ii) and (iii) go in analogy to that of Proposition 3 in \cite{Cwiszewski-Kokocki}.
(ii) corresponds also to results from \cite{Daners} (here we need a version with locally Lipschitz nonlinearity).
To see (ii), observe that if $x_n\to x_0$ in $E$ and $\lma_n\to \lma_0$, then putting
$u_n(s):=\Phi_s(x_n, \lma_n)$ for $s\in [0,T]$, one has
\begin{eqnarray*}\beta\left(\{ u_n(t)\}_{n\ge 1}\right)
\leq \beta(\{R^{(\lma_n)}(t,0)x_n\}_{n\ge 1}) + \beta\left( \left\{ \int_{0}^{t} h_{t,n} (s) \d s \ | \ n\ge 1 \right\}\right)
\end{eqnarray*}
where $h_{t,n} (s):=R^{(\lma_n)}(t,s)F(s,u_n(s),\lma_n)$ for
$s\in[0,t]$. In view of the $(F_2)_{par}$ and the Gronwall
inequality, there is $M_0\geq 0$ such that $\|u_n(s)\| \leq M_0$,
for all $s\in [0,T]$ and $n\geq 1$, which again by $(F2)_{par}$
gives $M_1\geq 0$ such that $\|h_{t,n} (s)\| \leq M_1$, for $s\in
[0,T]$ and $n\geq 1$. This allows applying Remark
\ref{07012009-1059} and  Lemmata \ref{19092008-0029} and
\ref{08012009-2310} and as a result, for any $t\in[0,T]$,
\begin{align*}
\beta\left(\{u_n(t)\}_{n\geq 1}\right) & \leq  \beta\left( \left\{ \int_{0}^{t} h_{t,n}(s) \d s \ | \ n\ge 1 \right\}\right) \\
& \leq \int_{0}^{t} \beta(\{h_{t,n} (s) \ | \ n\ge 1 \}) \d s
\leq k M \int_{0}^{t}\beta\left(\{u_n(s)\}_{n\ge 1}\right)
\end{align*}
By use of the Gronwall inequality, one gets $\beta\left(\{u_n(t)\}_{n\geq 1}\right)=0$ for any $t\in [0,T]$.
Hence, due to the fact that $u_n= \Sigma (x_n, w_{n}, \lma_n)$ with $w_n(s):=F(s,u_n (s),\lma_n)$
and Proposition \ref{30122008-1629} (ii) any subsequence of $(u_n)$ contains a subsequence
converging to some $u_0$. By Proposition \ref{30122008-1629} (i),
$u_0 = \Sigma (x_0, w_0, \lma_0)$ with $w_0(s):=F(s,u_0(s), \lma_0)$ and therefore
$u_0$ is a unique mild solution of (\ref{07012009-1523}). This completes the proof of (ii).\\
\indent (iii) Let $Q$ be an arbitrary bounded subset of $E$
and let $Q_0$ be a countable subset of $Q$ such that
$\overline{Q_0} \supset Q$ and $\Lambda_0$ a countable
dense subset of $[0,1]$. First observe that using $(F_2)_{par}$ and the Gronwall inequality as before,
we infer that the sets $\Phi_t (\overline{Q_0} \times [0,1])$, $t\in [0,T]$, are contained in a ball,
and clearly, by use of Lemmata \ref{08012009-2310} and
\ref{19092008-0029},
\begin{align*}
\beta(\Phi_t(Q_0\times \Lambda_0)) & \leq  \beta \left( \{ R^{(\lma)}(t,0)x\mid x\in Q_0, \, \lma\in \Lambda_0 \}\right) + \\
& \quad +\beta\left( \left\{ \int_{0}^{t}w_{x,\lma, t}(s) \d s \,\mid\,
x\in Q_0, \lma\in \Lambda_0 \right\}\right)\\
& \leq e^{-\omega t}\beta(Q_0) + \int_{0}^{t} \beta (\{w_{x,\lma,t}(s) \,\mid x\in Q_0, \lma\in \Lambda_0 \}) \d s,
\end{align*}
where $w_{x,\lma,t}(s):=R^{(\lma)}(t,s) F(s,\Phi_s(x,\lma),\lma)$
for $s\in [0,t]$. Further observe that, by (\ref{03082009-1958})
and Lemma \ref{08012009-2310}, for any $s\in [0,t]$, one has
\begin{align*}
\beta (\{w_{x,\lma,t}(s) \!\mid\! x\!\in\! Q_0, \lma\!\in\! \Lambda_0 \}) \\
& \hspace{-20mm}\leq \beta (\{R^{(\lma)}(t,s)z \, \mid \, z\!\in\! F([0,T]\!\times\! \Phi_s(Q_0\!\times\! \Lambda_0)\!\times\! \Lambda_0), \ \lma\!\in\!\Lambda_0\}) \\
& \hspace{-20mm}\leq e^{-\omega(t-s)} \beta\left( F([0,T]\!\times\! \Phi_s(Q_0\!\times\! \Lambda_0)\!\times\! \Lambda_0)\right)\\
& \hspace{-20mm}\leq k e^{-\omega(t-s)} \beta(\Phi_s(Q_0\!\times\! \Lambda_0)).
 \end{align*}
Combining the previous two inequalities together and applying the
Gronwall inequality give $\beta(\Phi_t(Q_0\times \Lambda_0)) \leq
e^{(k-\omega)t} \beta(Q_0)$ and finally $\beta(\Phi_t(Q\times
[0,1])) \le \beta(\Phi_t(\overline{Q_0}\times
\overline{\Lambda_0}))\leq \beta(\overline{\Phi_t(Q_0\times
\Lambda_0)}) =\beta(\Phi_t(Q_0\times \Lambda_0)) \leq
e^{(k-\omega)t} \beta(Q_0) \leq e^{(k-\omega)t} \beta(Q)$. \hfill
$\square$

Now we pass to the hyperbolic case. We shall assume in the rest of this section that
$V$ is a Banach space which is densely and continuously embedded into $E$.
Given a linear operator $A:D(A)\to E$ generating
a $C_0$ semigroup $\{S_A(t)\}_{t\geq 0}$ of bounded linear operators on $E$,
$V$ is said to be {\em $A$-admissible} provided $V$ is an invariant subspace for each $S_A(t)$ for $t\geq 0$ and the family of restrictions $\{S_{A}(t)_V:V\to V\}_{t\geq 0}$ ($S_A(t)_V x:=S_A(t)x$, $x\in V$)
is a $C_0$ semigroup on $V$. Define the \emph{part of $A$ in the space $V$} as  a linear operator $A_V:D(A_V)\to V$ given by $
D(A_V):=\{v\in D(A)\cap V\,\mid\, A v\in V\}$, $A_V v:=A v$ for $v\in D(A_V)$.
In view of \cite[Ch. 4, Theorem 5.5]{Pazy}, if $V$ is $A$-admissible then $A_V$ is the generator of the $C_0$ semigroup $\{S_A(t)_V\}_{t\geq 0}$.

\vspace{-5mm}

\begin{Prop}{\em (see \cite[Ch. 5, Theorem 3.1]{Pazy})}\label{17092008-1611}
Let $\{A(t)\}_{t\in [0,T]}$ be a family of linear operators on a Banach space $E$ satisfying the following conditions\\[2mm]
\noindent $(Hyp_1)$ \parbox[t]{140mm}{$\{A(t)\}_{t\in [0,T]}$ is
a stable family of infinitesimal generators of $C_0$ semigroups, i.e. there are $M\geq 1$ and $\omega\in\R$ such that
$$
\|S_{A(t_1)}(s_1)\ldots S_{A(t_n)}(s_n)\|_{{\cal L}(E,E)}\leq M e^{\omega(s_1+\ldots+s_n)},
$$
whenever $0\leq t_1 \leq \ldots \leq t_n\leq T$ and $s_1,\ldots, s_n\geq 0$, where $\{S_{A(t)}(s)\}_{s\geq 0}$ is the $C_0$ semigroup generated by $A(t)$;} \\[2mm]
\noindent $(Hyp_2)$ \parbox[t]{140mm}{$V$ is $A(t)$-admissible for each $t\in [0,T]$ and the family $\{ A_V (t)\}_{t\in [0,T]}$ is a stable family of generators of
$C_0$ semigroups with constants $M_V\geq 1$ and $\omega_V \in\R$;} \\[1mm]
\noindent $(Hyp_3)$ \parbox[t]{140mm}{$V\subset D(A(t))$ and $A(t)\in{\cal L}(V,E)$ for $t\in [0,T]$ and the mapping $[0,T]\ni t\mapsto A(t)\in {\cal L}(V,E)$ is continuous.}\\[1mm]
Then there exists a unique evolution system $\{R(t,s)\}_{0\leq s\leq t\leq T}$ in $E$ with the following properties\\[1mm]
\pari{(i)}{ $\|R(t,s)\| \leq M e^{\omega(t-s)}$ \quad for \ $0\leq s\leq t\leq T$;}\\[1mm]
\pari{(ii)}{ $\left.\frac{\partial^+}{\partial t} R(t,s)v\right|_{t=s}= A(s)v$
\quad for \ $v\in V$, $s\in [0,T)$;}\\[1mm]
\pari{(iii)}{ $\frac{\partial}{\partial s} R(t,s)v = -R(t,s)A(s)v$ \quad for \ $v\in V$, $0\leq s\leq t\leq T$.}
\end{Prop}

\vspace{-3mm}

Using homotopy invariants will require the continuity of linear
evolution systems with respect to parameters.

\vspace{-5mm}

\begin{Prop}\label{14012008-1812}
Let, for each $\lma\in [0,1]$, a family $\{A^{(\lma)}(t)\}_{0\leq t\leq T}$ satisfy conditions $(Hyp_1)-(Hyp_3)$ with constants $M, M_V, \omega, \omega_V$ independent of $\lambda$ and let $R^{(\lma)}=\{R^{(\lma)}(t,s)\}_{0\leq s\leq t\leq T}$ be the corresponding evolution systems in $E$ determined by Proposition \ref{17092008-1611}. If, for any $\lma_0\in [0,1]$,
\begin{equation}\label{13072009-1426}
\int_{0}^{T} \|A^{(\lma)}(\tau) - A^{(\lma_0)} (\tau)\|_{{\cal L}(V,E)} d \tau\to 0 \quad
\mbox{ as } \lma\to \lma_0,
\end{equation}
then $\{ R^{(\lma)}\}_{\lma\in [0,1]}$ is a continuous family of evolution systems in $E$.
\end{Prop}
{\bf Proof.} We use the construction from \cite[Ch. 5, Theorem 3.1]{Pazy}.
Recall that for any $\lma\in [0,1]$ and $x\in E$
$$
R^{(\lma)}(t,s)x := \lim_{n\to +\infty} R_{n}^{(\lma)}(t,s)x \quad\mbox{ for \ } 0\leq s\leq t\leq T,
$$
where, for each $n\geq 1$, the operator $R_{n}^{(\lma)} (t,s):E\to
E$ is given by (\footnote{Here we adopt the convention that
$\prod_{k=1}^n T_k := T_n\circ T_{n-1}\circ\ldots \circ T_1$, for
the sequence $T_1, T_2, \ldots, T_n$ of bounded operators on
$E$.})
$$
R_{n}^{(\lma)}(t,s)\!:= \!\!\left\{
\begin{array}{ll}
\!\! S^{(\lma)}_j(t\!-\!s) &  \mbox{ if }  s,t\in [t_{j}^{n}, t_{j+1}^{n}], \, s \le t,\\
\!\! S^{(\lma)}_k(t\!-\!t_{k}^{n})\!\left(\prod\limits_{j=l+1}^{k-1} \!\!S^{(\lma)}_{j}(T/n)\! \right)\!S^{(\lma)}_{l}(t_{l+1}^{n}\!-\!s)&
\mbox{ if } s\in [t_{l}^{n}, t_{l+1}^{n}], \, t\in [t_{k}^{n}, t_{k+1}^{n}],\\
&\mbox{ and } k>l \ge 0,
\end{array}
\right.
$$
with $t_{j}^{n}:=(j/n)T$, $S_{j}:=S_{A^{(\lma)}(t_{j}^n)}$, for $j=0,1,\ldots, n$.
Moreover recall (after \cite{Pazy}) that $\{R_{n}^{(\lma)}(t,s)\}_{0\leq s\leq t\leq T}$ are evolution systems such that
\begin{align}\label{14012009-1920}
\|R_{n}^{(\lma)}(t,s)\|_{{\cal L}(E,E)}\leq Me^{\omega(t-s)},
\quad \|R_{n}^{(\lma)}(t,s)\|_{{\cal L}(V,V)}\leq M_V e^{\omega_V
(t-s)}, \quad R_n^{(\lma)}(t,s)V\!\subset \! V,
\end{align}
for $0\leq s\leq t\leq T$ and for any $v\in V$
\begin{align}\label{14012009-1921}
\frac{\partial}{\partial t} R_{n}^{(\lma)}(t,s)v &
=A_{n}^{(\lma)}(t)R_n^{(\lma)}(t,s)v && \hspace{-10mm} \mbox{ for
\ } t\not\in \{t_{0}^{n},t_{1}^{n},\ldots, t_{n}^{n}\}, \ s\le t,
\\ \label{14012009-1922} \frac{\partial}{\partial s}
R_{n}^{(\lma)}(t,s)v & = - R_n^{(\lma)}(t,s) A_{n}^{(\lma)}(s)v
&&\hspace{-10mm}\mbox{ for \ }s\not\in
\{t_{0}^{n},t_{1}^{n},\ldots, t_{n}^{n}\}, \ s\le t,
\end{align}
with $A_{n}^{(\lma)}(t):=A^{(\lma)}(t_{k}^{n})$ if $t_k^n \leq t <
t_{k+1}^{n}$ for $k=0,\ldots,n-1$ and
$A_n^{(\lma)}(T):=A^{(\lma)}(T)$. Observe that in view of
$(Hyp_3)$ for any $\lma\in [0,1]$, one has $\|A_{n}^{(\lma)}(t) -
A^{(\lma)}(t)\|_{{\cal L}(V,E)} \to 0$ as  $n\to +\infty$
uniformly with respect to $t\in [0,T]$. Fix any $v\in V$, $\lma,
\mu\in [0,1]$, $n\geq 1$ and $s,t \in [0,T]$ with $s<t$ and define
$\phi:[s,t]\to E$ by $\phi(r):=
R_n^{(\lma)}(t,r)R_n^{(\mu)}(r,s)v$. In view of
(\ref{14012009-1920}), (\ref{14012009-1921}) and
(\ref{14012009-1922}) the map $\phi$ is differentiable on $[s,t]$
except finite number of points and
\begin{eqnarray*}
R_n^{(\mu)}(t,s)v- R_n^{(\lma)}(t,s)v = \phi(t)-\phi(s)
=\int_{s}^{t} \phi'(r) dr, \\
= \int_{s}^{t}  \left(R_n^{(\lma)}(t,r)(A_n^{(\mu)}(r))-A_n^{(\lma)}(r))R_n^{(\mu)}(r,s)v \right)dr.
\end{eqnarray*}
Hence, by (\ref{14012009-1920}),
$$
\| R_n^{(\mu)}(t,s)v- R_n^{(\lma)}(t,s)v \|
\leq MM_Ve^{(\omega +\omega_V)T}\|v\|_{V} \int_{0}^{T} \|A_n^{(\mu)}(r)-A_n^{(\lma)}(r)\|_{{\cal L}(V,E)} dr.
$$
Passing to the limit with $n\to +\infty$, we get
$$
\|R^{(\mu)} (t,s)v - R^{(\lma)}(t,s)v\| \leq
MM_Ve^{(\omega +\omega_V)T}\|v\|_{V} \int_{0}^{T} \|A^{(\mu)}(r)-A^{(\lma)}(r)\|_{{\cal L}(V,E)} dr
$$
and in consequence, $R^{(\lma)}(t,s)v\to R^{(\lma_0)}(t,s)v$ for any $v\in V$, as $\lma\to \lma_0$ and the convergence is uniform with respect to $s,t$.
Using the density of $V$ in $E$, we complete the proof since $\|R^{(\lambda)}(t,s)\| \leq M e^{\omega(t - s)}$ for $\lambda\in[0,1]$ and $0\leq s\leq t\leq T$.  \hfill $\square$

The following  criterion for verification of conditions $(Hyp_1)-(Hyp_3)$ is useful in applications.
\begin{Prop} {\em (\cite[Ch. 5, Theorem 4.8]{Pazy})}
\label{27022009-1214} Suppose that a family $\{A(t) \}_{t\in
[0,T]}$, where $D(A(t))=D$ for any $t\in [0,T]$ and some $D\subset
E$, is stable and for each $v\in D$ the mapping $[0,T]\ni t\mapsto
A(t)v\in E$ is continuously differentiable. Then, the family $\{
A(t)\}_{t\in [0,T]}$ satisfies conditions $(Hyp_1)-(Hyp_3)$ with
$V:=D$ equipped with the norm given by $\|v\|_V:=\|A(0)v\| +
\|v\|$ for $v\in V$.
\end{Prop}

\section{Averaging method for periodic solutions}

We shall deal with the periodic problem
$$
\left\{\begin{array}{ll}
\dot u(t) = A(t)u(t) + F(t,u(t)), \quad t\in [0,T] \\
u(0) = u(T)
\end{array}\right.\leqno{(P)}
$$
where the family $\{A(t)\}_{t\in [0,T]}$ of linear operators  on a separable Banach space $E$ satisfies
a more restrictive variant of $(Hyp_1)$ (from Proposition \ref{17092008-1611})\\ [1mm]
\noindent $(Hyp'_1)$ \ \parbox[t]{140mm}{there is $\omega > 0$, such that
$$
\|S_{A(t_1)}(s_1)\ldots S_{A(t_n)}(s_n)\|_{{\cal L}(E,E)}\leq
e^{-\omega(s_1+\ldots+s_n)};
$$
whenever $0\leq t_1 \leq \ldots \leq t_n\leq T$ and $s_1,\ldots, s_n\geq 0$,} \\[2mm]
conditions $(Hyp_2)$, $(Hyp_3)$ and, additionally, \\[1mm]
$(Hyp_4)$ \parbox[t]{135mm}{ there is $\mu_0 > -\omega$ such that the space $(\mu_0 I - A_0)V$ is dense in $E$,
where $$A_0:= \frac{1}{T}\int_{0}^{T} A(\tau)\, d\tau \in {\cal L}(V,E);$$}\\
$(Hyp_5)$ \parbox[t]{135mm}{ $A(0)x = A(T)x$ for $x\in D(A(0)) = D(A(T))$.}\\[2mm]
\noindent Furthermore, we assume that a continuous mapping $F:[0,T] \times E\to E$\\
\noindent $(F_{1})$ \parbox[t]{138mm}{is locally Lipschitz with respect to the second variable uniformly with respect to the first one;}\\
\noindent $(F_{2})$ \parbox[t]{138mm}{ has sublinear growth uniformly with respect to
the first variable, i.e. there is a constant $c>0$ such that
$$\|F(t,x)\|\le c(1 + \|x\|)\quad\mbox{ for $x\in E$, $t\in [0,T]$};$$}\\
\noindent $(F_{3})$ \parbox[t]{138mm}{ there is $k\in[0,\omega)$
such that
$$\beta(F([0,T]\times Q))\le k\beta(Q) \quad\mbox{ for any bounded
$Q\subset E$};$$}\\
\noindent $(F_{4})$ \parbox[t]{138mm}{  $F(0,x) = F(T,x)$ for $x\in E$.}\\[2mm]
The existence of periodic solutions will be obtained by means of a
continuation principle for a parameterized family of periodic
problems
$$
\left\{\begin{array}{ll}
\dot u(t) = \lambda A(t)u(t) + \lambda F(t,u(t)),\quad t\in [0,T] \\
u(0) = u(T)
\end{array}\right.\leqno{(P_{T,\lambda})}
$$
with the parameter $\lambda\geq 0$.
For any $x\in E$ and $\lma\in [0,T]$, by $u(\,\cdot\,;x,\lambda)$ denote the unique mild solution of
\begin{eqnarray}\label{29082008-2340}
\dot u(t)=\lambda A(t)u(t)+\lambda F(t,u(t)), \quad t\in [0,T]
\end{eqnarray}
satisfying the initial condition $u(0;x,\lma) = x$. The translation along trajectories operator
for $(\ref{29082008-2340})$ is denoted by $\Phi_{t}^{(\lambda)}:E \to E$ where $t\in [0,T]$.
A point $(x,\lambda)\in E\times [0,+\infty)$ is a {\em $T$--periodic point}
for (\ref{29082008-2340}) if $\Phi_{T}^{(\lambda)}(x)=x$.
We say that $x_0\in E$ is  a {\em branching point } (or  a {\em cobifurcation point}) for $(P_{T,\lambda})$, $\lma\geq 0$,
if there exists a sequence of $T$-periodic points $(x_n,\lambda_n)\in E \times(0,+\infty)$
such that $\lambda_n\to 0$ and $x_n\to x_0$ as $n\to +\infty$.

\begin{Th}\label{18092008-1450}
If $x_0\in E$ is a branching point of $(P_{T,\lambda})$, then
$\widehat{A} x_0 + \widehat{F} (x_0) = 0$
where
$\widehat A:D(\widehat A) \to E$ is the closure of $A_0$
and $\widehat F:E\to E$ is given by $\widehat F(x):=(1/T)\int_{0}^{T} F(\tau,x)\,d\tau$.
\end{Th}
\begin{Rem}\label{19072009-2137}{\em
(a) Recall that due to \cite[Ch. 1, Th. 4.3]{Pazy}, if $A:D(A)\to E$ generates a $C_0$ semigroup of contractions,
then the dissipativity condition
\begin{equation}\label{rr1}
\|x-\lma Ax\| \geq \|x\| \quad \mbox{ for any \ } x\in D(A), \ \lambda > 0
\end{equation}
is equivalent to
$$
(p,Ax)\leq 0 \quad\mbox{ for any \ } x\in D(A), \ p\in J(x),
$$
where $J(x):=\{ p \in E^* \mid \la p, x\ra = \|x\|^2 = \|p\|^{2}\}$ is the dual set of $x$.\\
\indent (b) Hence, if $(Hyp'_1)$ and $(Hyp_3)$ hold, then $\omega
I + A_0= \omega I + (1/T)\int_{0}^{T} A(\tau) d\tau$ in ${\cal
L}(V,E)$ is a dissipative operator. This implies that the closure
$\widehat A:D(\widehat A) \to E$ of $A_0$ is a well-defined linear
operator and, by (\ref{rr1}), the operator $\widehat
A_\omega:=\omega I + \widehat A$ is also dissipative, hence
$\lambda I - \widehat A_\omega$ has closed range whenever $\lambda
> 0$. If condition $(Hyp_4)$ holds, then the operator $(\mu_0 +
\omega)I - \widehat A_\omega = \mu_0 I - \widehat A$ has closed
and dense range, since $\mu_0 + \omega > 0$.
It means that $(\mu_0 + \omega)I - \widehat A_\omega$ is $m$--dissipative, since its range is the whole $E$ and, by the Lumer-Phillips theorem, $\widehat A$ generates a $C_0$ semigroup such that $\|S_{\widehat A}(t)\|\le e^{-\omega t}$ for $t\geq 0$.\\
\indent (c) In particular, $(-\omega,+\infty)\subset\varrho(\widehat A)$ and, for each $\mu > -\omega$,
$$
(\mu I - A_0)V = (\mu I - \widehat A)V = (\mu I - \widehat A)(\mu_0 I - \widehat A)^{-1}V_0
$$
with $V_0 := (\mu_0 I - \widehat A)V$ being a dense subset of $E$.
Since $(\mu I - \widehat A)(\mu_0 I - \widehat A)^{-1}:E\to E$
is a bounded bijection, we infer that $(\mu I - A_0)V$ is dense in $E$ for any $\mu > -\omega$.\\
\indent (d) If $\{A(t)\}_{t\in [0,T]}$ satisfies conditions $(Hyp'_1)$, $(Hyp_2)$ and $(Hyp_3)$, then, in view of Proposition \ref{17092008-1611}
and point (b), one has
$$\|R(t,s)\|_{{\cal L}(E,E)}\leq e^{-\omega(t-s)}
\quad \mbox{ for any } s,t\in [0,T], s\leq t.
$$}
\end{Rem}
In the proof of Theorem \ref{18092008-1450} and later on in the section we use the following lemma.
\begin{Lem}\label{14092008-1339}
Let $\{A(t)\}_{t\in [0,T]}$ satisfy $(Hyp'_1)$, $(Hyp_2)$ -- $(Hyp_4)$ and let, for each $\mu\in [0,1]$, the family $\{A^{(\mu)}(t)\}_{t\in [0,T]}$ be defined by $A^{(\mu)}(t):=-\mu I + (1-\mu)A(t)$ for $t\in [0,T]$.
Then, for each $\lma\geq 0$ and $\mu\in [0,1]$, the family of operators $\{\lma A^{(\mu)}(t)\}_{t\in [0,T]}$
satisfies $(Hyp_1)'$, $(Hyp_2)$ -- $(Hyp_4)$ and the corresponding evolution systems $\{R^{(\mu,\lambda)}(t,s)\}_{0\leq s\leq t\leq T}$, $\lambda\ge 0$, $\mu\in [0,1]$ have the following properties \\
\pari{(i)}{for any $x\in E$, $t,s\in [0,T]$ with $s\leq t$,
$(\lambda_n)$ in $(0,+\infty)$ and $(\mu_n)$ in $[0,1]$ such that $\lma_n\to 0$, $\mu_n \to \mu_0$, one has
$$R^{(\mu_n,\lambda_n)}(t,s)x \to x \quad \mbox{ as } n\to +\infty,$$
uniformly with respect to $t,s\in [0,T]$ with $s \le t$;}\\[2mm]
\pari{(ii)}{if $(k_n)$ is a sequence of positive integers and sequences $(\lma_n)$ in $(0,+\infty)$ and $(\mu_n)$ in $[0,1]$ are such that $k_n\to +\infty$, $k_n\lambda_n\to \eps$ for some $\eps>0$ and $\mu_n\to \mu_0$ for some $\mu_0\in [0,1]$, then for any $x\in E$
$$
R^{(\mu_n,\lambda_n)}(T,0)^{k_n} x\to S_{\widehat{A^{(\mu_0)}}}(\eps T)x \quad \mbox{ as } n\to\infty,
$$
where $\widehat{A^{(\mu_0)}}$ is the closure of the operator
$\frac{1}{T}\int_{0}^{T} A^{(\mu_0)}(\tau)\,d\tau$;}\\[2mm]
\pari{(iii)}{if $(k_n)$, $(\lma_n)$ and $(\mu_n)$ are as in {\em (ii)}, then for any $x\in E$
$$\lambda_n(I+R^{(\mu_n,\lambda_n)}(T,0)+\ldots+R^{(\mu_n,\lambda_n)}(T,0)^{k_n-1})x\to
\frac{1}{T}\int\limits_{0}^{\eps T}\!\! S_{\widehat{A^{(\mu_0)}}}(\tau)x \,d\tau \quad \mbox{ as } n\to+\infty.$$} \\[-10mm]
\end{Lem}
\noindent\textbf{Proof.} (i) It is easy to check that, for each
$\lma > 0$ and $\mu\in [0,1]$, the family $\{\lma
A^{(\mu)}(t)\}_{t\in [0,T]}$ satisfies $(Hyp'_1)$ with constant
$\omega:=\lambda \min\{1,\omega\}$ and conditions $(Hyp_2)$ --
$(Hyp_3)$ as well. From now on we write $A_{0}^{(\mu)} := -\mu I +
(1-\mu) A_0$ for $\mu\in [0,1]$. We claim that also $(Hyp_4)$
holds. Indeed if $\mu = 1$, then $A^{(\mu)}(t) = I$ for $t\in
[0,T]$ and $(a_{\mu, \lambda} I - A^{(\mu)}_0)V$ with $a_{\mu,
\lambda} = 0$ is dense in $E$, for $\lambda > 0$. If $\mu \neq 1$,
then putting $a_{\mu, \lambda} := (1 - \mu)\lambda\mu_0 -
\lambda\mu$ we see that $a_{\mu, \lambda} > - \omega$, since
$\mu_0 > -\omega$, and $(a_{\mu, \lambda} I - \lambda
A^{(\mu)}_0)V = \lambda(1 - \mu)(\mu_0 I - A_0)V$ is dense in $E$
and thus $\lambda A^{(\mu)}_0$ satisfies $(Hyp_4)$. For the
corresponding evolution system $R^{(\mu,\lma)}$, one gets, for any
$(\lma,\mu)\in (0,+\infty)\times [0,1]$,
\begin{equation}\label{18072009-2214}
\|R^{(\mu,\lambda)}(t,s)\|\leq e^{-\lambda\overline{\omega} (t - s)} \leq 1 \quad\mbox{ for \ }t,s\in [0,T], \, s\le t,
\end{equation}
with $\overline\omega:=\min\{1,\omega\}$, and
$\frac{\partial}{\partial s} R^{(\mu,\lma)}(t,s)v = -\lma R^{(\mu,\lma)}(t,s)A^{(\mu)}(s)v$, for $v\in V$, $0\leq s\leq t\leq T$.
In consequence, for any $v\in V$,  $t,r\in [0,T]$, $r \le t$, $\mu\in [0,1]$ and $\lma>0$, one has
$$
R^{(\mu,\lma)}(t,r)v-v= R^{(\mu,\lma)}(t,r)v -
R^{(\mu,\lma)}(t,t)v=\lambda\int_{r}^{t}
R^{(\mu,\lma)}(t,s)A^{(\mu)}(s)v \d s.
$$
Since, for any $s\in [r,t]$,
$$
\|R^{(\mu,\lma)}(t,s)A^{(\mu)}(s)v\| \leq \|R^{(\mu,\lma)}(t,s)\|\|A^{(\mu)}(s)\|_{{\cal L}(V,E)} \|v\|_{V}
 \leq \|A^{(\mu)}(s)\|_{{\cal L}(V,E)} \|v\|_{V},
$$
we infer that $\|R^{(\mu,\lma)}(t,r)v-v\| \leq \lma C\|v\|_{V}$
with $C:=\sup_{\mu\in[0,1]} \int_{0}^{T} \|A^{(\mu)}(s)\|_{{\cal
L}(V,E)} \d s<+\infty$. This, due to the density of $V$ in $E$,
means that
$$
\lim_{\lma\to 0^+,\,  \mu\to \mu_0 } R^{(\mu,\lma)}(t,r)x = x \quad\mbox{ for any } x\in E
$$
uniformly with respect to $t,r\in[0,T]$ with $r\le t$,  which implies (i).\\
\indent Define a map $L:(0,+\infty)\times [0,1]\to {\cal L}(E,E)$ by
$$
L(\lambda,\mu):= R^{(\mu,\lambda)}(T,0) \quad \mbox{ for any } \lambda>0 \mbox{ and }\mu\in [0,1].
$$
Clearly, by (\ref{18072009-2214}), $\|L(\lambda,\mu)\| \leq 1$ for each $(\lma,\mu)\in (0,+\infty)\times [0,1]$.
Observe also that, for each $\mu\in[0,1]$, $a_\mu := a_{\mu, \lambda} /\lambda$ is such that $(a_\mu I - A_{0}^{(\mu)})V$ is dense in $E$. Further, for any $v\in V$,
\begin{eqnarray}
\lma^{-1}(L(\lma,\mu)v-v)=\lma^{-1}(R^{(\mu,\lma)}(T,0)v-v) =
\int_{0}^{T} R^{(\mu,\lma)}(T,s) A^{(\mu)}(s)v \d s
\end{eqnarray}
and
\begin{align*}
&\| \lma^{-1}(L(\lma,\mu)v-v) - T A^{(\mu_0)}_0v\|\leq
\int_{0}^{T} \|R^{(\mu,\lma)}(T,s)A^{(\mu)}(s)v-A^{(\mu_0)}(s)v\| \d s\\
&\quad \leq \int_{0}^{T} \|R^{(\mu,\lma)}(T,s) A^{(\mu)}(s)v-A^{(\mu)}(s)v\| \d s +
\|v\|_{V}\int_{0}^{T} \|A^{(\mu)} (s)-A^{(\mu_0)}(s)\|_{{\cal L}(V,E)} \d s,
\end{align*}
which by use of point (i) of this lemma and $(Hyp_3)$ gives
$$
\lim_{\lma\to 0^+,\, \mu\to \mu_0} \lma^{-1}(L(\lma,\mu)v-v) = TA^{(\mu_0)}_0v \quad \mbox{ for \ } v\in V.
$$
Hence, applying Theorem \ref{29082008-1821} and changing the time variable we get (ii)
and (iii) as the closure of $A^{(\mu_0)}_0$ is equal to $\widehat {A^{(\mu_0)}}$. \hfill $\square$

\noindent {\bf Proof of Theorem \ref{18092008-1450}.}
Let sequences $(\lambda_n)$ in $(0,+\infty)$ and $(x_n)$ in $E$ be
such that $\lambda_n \to0$, $x_n\to x_0$ as $n\to \infty$ and $\Phi_T^{(\lambda_n)}(x_n)=x_n$ for each $n\geq 1$. Then, by definition
\begin{equation}\label{14092008-1236}
\Phi_t^{(\lambda_n)} (x_n) = R^{(\lambda_n)}(t,0)x_n + \lambda_n
\int_{0}^{t} R^{(\lambda_n)}(t,s)F(s,\Phi_s^{(\lambda_n)}(x_n)) \d
s
\end{equation}
for any $t\in [0,T]$, where $\{R^{(\lambda)}(t,s)\}_{0\leq s\leq t\leq T}$ denotes the evolution system generated by the family $\{\lambda A(t)\}_{t\in [0,T]}$.
This yields
\begin{align*}
\|x_n-\Phi_t^{(\lambda_n)}(x_n)\| & = \|R^{(\lambda_n)}(t,0)x_n - x_n\| +
\left\|\lambda_n \int_{0}^{t} R^{(\lambda_n)}(t,s)F(s,\Phi_s^{(\lambda_n)}(x_n)) \d s \right\|\\
& \le \|R^{(\lambda_n)}(t,0)x_n - x_n\|+ \lambda_nc \int_0^t (1 +
\|\Phi_s^{(\lambda_n)}(x_n)\|)\d s.
\end{align*}
In view of Lemma \ref{14092008-1339} (i) with $\mu_n := 0$ for $n\ge 1$ and the boundedness of
$\{\Phi_s^{(\lambda_n)}(x_n)\mid s\in [0,T], n\geq 1 \}$, we infer that
$\Phi_t^{(\lambda_n)}(x_n)\to x_0$, uniformly with respect to $t\in [0, T]$.
Further, by (\ref{14092008-1236}), one has
\begin{equation}
x_n = \Phi_T^{(\lambda_n)}(x_n) = R^{(\lambda_n)}(T,0)x_n +
\lambda_n\int_{0}^{T}R^{(\lambda_n)}(T,s)F(s,\Phi_s^{(\lambda_n)}(x_n))
\d s
\end{equation}
and consequently, for each $k\geq 0$
$$
R^{(\lambda_n)}(T,0)^{k} x_n = R^{(\lambda_n)}(T,0)^{k+1} x_n +
\lambda_n R^{(\lambda_n)}(T,0)^{k}\int_0^T R^{(\lambda_n)}(T,s)F(s,\Phi_s^{(\lambda_n)}(x_n))
d s.
$$
Let $\eps>0$ be arbitrary and let $(k_n)$ be a sequence
of positive integers such that $k_n\lambda_n\to \eps$.
Summing up the above equalities with $k=0,1,\ldots, k_n-1$ for any $n\geq 1$, we obtain
$$
\begin{array}{l}
x_n = R^{(\lambda_n)}(T,0)^{k_n}x_n + \lambda_n\left[\sum\limits_{k=0}^{k_n-1} R^{(\lambda_n)}(T,0)^{k}\right]
\left(\int\limits_{0}^{T} R^{(\lambda_n)}(T,s)F(s,\Phi_s^{\lambda_n} (x_n) )
\d s\right)
\end{array}
$$
and, by use of Lemma \ref{14092008-1339} (ii) and (iii) with $\mu_n := 0$ for $n\ge 1$,
$$x_0 = S_{\widehat{A}} (\eps T) x_0 + \left[\frac{1}{T}\int_0^{\eps T} S_{\widehat{A}}(\tau) d\tau \right] \left(\int_{0}^{T} F(s,x_0)\d s\right).$$
In consequence
$$
-\frac{1}{\eps T} \left( S_{\widehat{A}} (\eps T) x_0 - x_0\right) =
\frac{1}{\eps T}\int_0^{\eps T}
S_{\widehat{A}}(\tau)\widehat{F} (x_0) d\tau.
$$
Thus, since $\eps\!>\!0$ was arbitrary, letting $\eps\!\to\! 0^+$, one has
$-\widehat A x_0 \!=\! \widehat F(x_0)$. \hfill $\square$

\begin{Th}\label{18092008-2220} {\em (Averaging principle)}
Let $\{A(t)\}_{t\in [0,T]}$ be a family of generators of $C_0$ semigroups satisfying
$(Hyp'_1)$, $(Hyp_2)$ -- $(Hyp_5)$ and let $F:[0,T] \times E \to E$
be a continuous map with properties $(F_1)$ -- $(F_4)$.
If $U\subset E$ is an open bounded set such that $\widehat{A}x + \widehat{F}(x)\neq 0$ for any
$x \in\partial U\cap D(\widehat{A})$, then there exists $\lambda_0>0$
such that for all $\lambda\in (0,\lambda_0]$, $\Phi_T^{(\lambda)} (x)\neq x$
for all $x\in\partial U$ and
\begin{equation}
         \Deg (\widehat{A}+\widehat{F}, U) = \deg (I-\Phi_T^{(\lambda)}, U).
\end{equation}
Here $\deg$ stands for the topological degree for condensing vector fields (see {\em \cite{Akhmerov-etal}} or {\em \cite{Nussbaum}})
and $\Deg(\widehat{A}+\widehat{F}, U):= \deg(I+{\widehat A}^{-1} \widehat F, U)$ (see {\em \cite{Cwiszewski-Kokocki}}).
\end{Th}
In the proof we shall need the following lemmata.
\begin{Lem}\label{19092008-0012}{\em (see \cite[Lemma 5.4]{Cwiszewski-Kokocki})}
Let $T_n:E\to E$, $n\geq 1$, be bounded linear operators, such that, for any $x\in E$,
$(T_n x)$ is a Cauchy sequence {\em(}\footnote{This is actually equivalent to the existence of a bounded operator
$T:E\to E$ being the strong limit of $(T_n)$.}{\em)}. Then, for any bounded set $\{x_n\}_{n\geq 1} \subset E$
$$
\beta\left( \{T_n x_n\}_{n\geq 1} \right) \leq \left( \limsup_{n\to +\infty}
\|T_n\|\right) \beta\left( \{x_n\}_{n\geq 1}\right).
$$
\end{Lem}
\begin{Lem}{\em(cf. Step 2 in the proof  of Theorem 5.1 in \cite{Cwiszewski-Kokocki})} \label{22092008-1201}
Let $A$ be a generator of a $C_0$ semigroup $S_A$ such that $\|S_A(t)\|\leq e^{-\omega t}$ for $t\geq 0$
and $F:E \to E$ be a continuous map with $k\in [0,\omega)$
such that $\beta(F(Q))\leq k\beta(Q)$ for any bounded $Q$.
If an open  bounded $U\subset E$ is such that $Ax+F(x)\neq 0$ for each $x\in \partial U\cap D(A)$,
then there exists a locally Lipschitz compact mapping $F_L:E\to E$ such that
\begin{equation}\label{19072009-2232}
Ax+(1-\mu)F(x) + \mu F_L(x)\neq 0  \quad \mbox{ for \ } x\in \partial U\cap D(A), \ \mu\in [0,1].
\end{equation}
\end{Lem}
\begin{Lem}\label{19092008-0014}
Let $\{A^{(\mu)}(t)\}_{t\in [0,T]}$ for $\mu\in [0,1]$, satisfy
$(Hyp'_1)$, $(Hyp_2)-(Hyp_5)$ with the common, independent of
$\mu$, constants $\omega > 0$, $\omega_V$, $M_V$ and let
$F:[0,T]\times E \times [0,1]\to E$ be a continuous mapping
satisfying $(F_1)_{par}$ -- $(F_3)_{par}$ and the periodicity
condition
$$F(0,x,\lma)=F(T,x,\lma) \quad \mbox{ for \ }(x,\lma) \in E
\times [0,1].$$ Suppose that $U\subset E$ is open bounded and
$\widehat{A^{(\mu)}}x+\widehat{F}(x,\mu)\neq 0$ for $x\in \partial
U \cap D(\widehat{A^{(\mu)}})$ and $\mu\in [0,1]$, where
$\widehat{A^{(\mu)}}$ is the closure of $(1/T)\int_{0}^{T}
A^{(\mu)}(s) \d s$ and $\widehat{F}:E\times [0,1]\to E$ is given
by $\widehat{F}(x,\mu):=(1/T)\int_{0}^{T} F(s,x,\mu) \d s$. Then,
there exists $\lma_0>0$, such that, for any $\lma\in (0,\lma_0]$,
$$
\Psi_T^{(\lma)}(x,\mu)\neq x \quad \mbox{ for all } x\in \partial U, \ \mu\in [0,1],
$$
where $\Psi_T^{(\lma)}: \overline U\times [0,1]\to E$ is given by
$\Psi_{T}^{(\lma)}(x,\mu):=u(T;x,\mu,\lma)$ for $(x,\mu)\in \overline U\times [0,1]$, $\lma>0$
and $u(\,\cdot \,;x,\mu,\lma):[0,T]\to E$ is the unique mild solution of
$$\left\{
\begin{array}{l}\dot u(t)=\lma A^{(\mu)}(t) u(t)+ \lma F(t,u(t),\mu), \quad t\in [0,T] \\
u(0)=x.\end{array}\right.$$
\end{Lem}
{\bf Proof. }Suppose to the contrary that there exist sequences $(\lma_n)$ in $(0,+\infty)$ with $\lma_n \to 0^+$, $(x_n)$ in $\partial U$ and $(\mu_n)$ in $[0,1]$ such that $\Psi_{T}^{(\lma_n)}(x_n,\mu_n)=x_n$ for $n\geq 1$.
Without loss of generality, we may assume that $\mu_n\to \mu_0$ as
$n\to +\infty$, for some $\mu_0\in [0,1]$. Let $\{\overline
A^{(\mu)}(t)\}_{t\in [0,2T]}$, $\mu\in [0,1]$ and a mapping $\overline F: [0,2T]\times E\times [0,1]\to E$ be given by
\begin{align*}
\overline  A^{(\mu)}(t) & :=A^{(\mu)}(t-[t/T]T) && \hspace{-15mm}\mbox{ for \ } (t,\mu)\in [0,2T]\times [0,1], \\
\overline F (t,x,\mu) & :=F(t-[t/T]T,x,\mu) && \hspace{-15mm}\mbox{ for \ } (t,x,\mu)\in [0,2T]\times E\times [0,1].
\end{align*}
where $[s]$ stands for the integer part of $s\in\R$. It is easy to
check that, for each $\lambda\in(0,\infty)$ and $\mu\in [0,1]$,
the family $\{\lambda \overline A^{(\mu)}(t)\}_{t\in [0,2T]}$ and
the mapping $\lambda F$ satisfies $(Hyp'_1)$, $(Hyp_2)$ --
$(Hyp_4)$ and $(F1)_{par}$ -- $(F3)_{par}$. Denote by $\{\overline
R^{(\mu,\lma)}\}_{0\le s\le t\le 2T}$ the corresponding evolution
system obtained by Proposition \ref{17092008-1611}. From the very
construction of hyperbolic evolution systems (see \cite[Ch. 5,
Theorem 3.1]{Pazy} and the proof of Proposition
\ref{14012008-1812}), we see that, for all $\lma > 0$, and $\mu\in
[0,1]$,
\begin{equation}\label{19072009-1921}
\overline R^{(\mu,\lma)}(T+t,T+s) = \overline R^{(\mu,\lma)}(t,s) = R^{(\mu,\lma)}(t,s) \quad\mbox{ for \ } t,s\in[0,T], \ s\le t.
\end{equation}
For each $t\in [0,2T]$, define $\overline{\Psi}_{t}^{(\lma)}:E\times
[0,1]\to E$, by
$\overline{\Psi}_{t}^{(\lma)}(x,\mu):=\overline u(t; x,\mu,\lma)$,
where $\overline u(\, \cdot \,;x,\mu,\lma):[0,2T]\to E$ is a solution of
$$\left\{\begin{array}{l}
\dot u(t)=\lma\overline A^{(\mu)}(t)u(t)+\lma\overline F(t,u(t),\mu),\quad t\in [0,2T] \\
u(0) = x.
\end{array}\right.$$
It clearly follows from (\ref{19072009-1921}) that, for any $n\geq 1$ and $t\in [0,T]$,
$$
\Psi_{t}^{(\lma_n)} (x_n, \mu_n)=
\overline{\Psi}_{t+T}^{(\lma_n)}(x_n, \mu_n ) = \overline
R^{(\mu_n, \lma_n)}(T+t,t)\Psi_{t}^{(\lma_n)}(x_n,\mu_n) +\lma_n
\int_{t}^{T+t} w_{n,t}(s) \d s,
$$
with $w_{n,t}(s):=\overline R^{(\mu_n,\lma_n)}(T+t,s) \overline F(s,\overline{\Psi}_{s}^{(\lma_n)}(x_n,\mu_n),\mu_n)$
for $s\in [t,T+t]$. Therefore, for any integer $k\geq 0$, one has
\begin{align*}
\overline R^{(\mu_n,\lma_n)}(T+t,t)^{k}\Psi_{t}^{(\lma_n)}(x_n,\mu_n) & = \overline R^{(\mu_n,\lma_n)}(T+t,t)^{k+1}
\Psi_{t}^{(\lma_n)}(x_n,\mu_n)\\
& \qquad + \lma_n \overline R^{(\mu_n, \lma_n)}(T+t,t)^{k}
\int_{t}^{T+t} w_{n,t}(s) \d s.
\end{align*}
Putting $k_n:=[1/\lma_n]$ and summing up the above equalities with $k=0,\ldots, k_n-1$,
we find that
\begin{equation}\label{19072009-1956}
\Psi_{t}^{(\lma_n)}(x_n,\mu_n) = \overline
R^{(\mu_n,\lma_n)}(T+t,t)^{k_n}\Psi_{t}^{(\lma_n)}(x_n,\mu_n) +
K_n \left(\int_{t}^{T+t} w_{n,t}(s) \d s \right)
\end{equation}
where
$$
K_n:=\lma_n \sum_{k=0}^{k_n - 1}\overline
R^{(\mu_n,\lma_n)}(T+t,t)^{k} \quad\mbox{ for \ }n\geq 1.$$ By
(\ref{19072009-1921}) and the fact that $\lambda_n k_n \to 1$ as
$n\to +\infty$, going along the lines of the proof of Lemma
\ref{14092008-1339}, we infer that, for any $x\in E$,
$$
\overline R^{(\mu_n,\lma_n)}(T+t,t)^{k_n}x\to
S_{\widehat{A^{(\mu_0)}}}(T) x \quad \mbox{ and } \quad K_n x \to
\frac{1}{T}\int_{0}^{T}  S_{\widehat{A^{(\mu_0)}}}(s)x \d s.
$$
This, along with Lemmata \ref{19092008-0012} and \ref{19092008-0029},
gives
\begin{align*}
\beta\left(\left\{\Psi_{t}^{(\lma_n)}(x_n,\mu_n)\right\}_{n\ge 1}\right) \le \\
& \hspace{-30mm}\leq e^{-\omega T} \beta\left(\left\{\Psi_{t}^{(\lma_n)}(x_n,\mu_n)\right\}_{n\ge 1}\right)
+ \frac{1-e^{-\omega T}}{\omega T} \beta\left(\left\{\int_{t}^{T+t} w_{n,t}(s) \d s \right\}_{n\ge 1}\right) \\
& \hspace{-30mm}\leq e^{-\omega T} \beta\left(\left\{\Psi_{t}^{(\lma_n)}(x_n,\mu_n)\right\}_{n\ge 1}\right) +
\frac{1-e^{-\omega T}}{\omega T} \int_{t}^{T+t} \beta(\{w_{n,t}(s)\}_{n\ge 1}) \d s\\
& \hspace{-30mm} \leq e^{-\omega T} \beta\left(\left\{\Psi_{t}^{(\lma_n)}(x_n,\mu_n)\right\}_{n\ge 1}\right) +
\frac{1-e^{-\omega T}}{\omega T} \int_{t}^{T+t} k\beta\left(\left\{\overline{\Psi}_{s}^{(\lma_n)}(x_n,\mu_n)\right\}_{n\ge 1}\right)ds
\end{align*}
and, in consequence,
\begin{equation}
\beta\left(\left\{\Psi_{t}^{(\lma_n)}(x_n,\mu_n)\right\}_{n\ge 1}\right) \leq \frac{k}{\omega T}\int_{t}^{T+t} \beta\left(\left\{\overline{\Psi}_{s}^{(\lma_n)}(x_n,\mu_n)\right\}_{n\ge 1}\right) ds.
\end{equation}
Define $\phi:[0,2T]\to \R$ by
$$
\phi(s):= \beta\left(\left\{\overline{\Psi}_{s}^{(\lma_n)}(x_n,\mu_n)\right\}_{n\geq 1}\right) \quad \mbox{ for \ } s\in [0,2T].
$$
We claim that $\phi \equiv 0$. Indeed, otherwise  $M:=\sup_{s\in [0,2T]} \phi(s) \in (0,+\infty)$ and by its $T$-periodicity,  for $\eps\in (0,(1-k/\omega)M)$ there exists $t_\eps\in [0,T]$ such that
$$
M-\eps <\phi(t_\eps) \leq \frac{k}{\omega T}
\int_{t_\eps}^{T+t_\eps}\phi(s) \d s \leq (k/\omega) M<M-\eps,
$$
which is a contradiction. Hence, in particular $\beta(\{x_n \}_{n\ge 1})=0$ and without loss of generality we may assume that $x_n\to x_0$ as $n\to +\infty$, for some $x_0\in \partial U$.\\
\indent Further, observe that $\Psi_{t}^{(\lma_n)}(x_n,\mu_n) \to x_0$ as $n\to +\infty$, uniformly with respect to $t\in [0,T]$, which follows from the inequality
$$
\|\Psi_{t}^{(\lma_n)}\! (x_n,\mu_n)-x_0\|\!\leq\!
\|R^{(\mu_n,\lma_n)}\! (t,0) x_n \!- \!x_0\| \!+ \!\lma_n \!\int_{0}^{t}\!\!\! \|R^{(\mu_n,\lma_n)}(t,s)F(s,\!\Psi_{s}^{(\lma_n)}\!(\!x_n, \mu_n),\!\mu_n) \| \d s
$$
and Lemma \ref{14092008-1339} (i). Note that, for any $k\geq 0$
\begin{eqnarray*}
 R^{(\mu_n,\lma_n)}(T,0)^{k}x_n  = R^{(\mu_n,\lma_n)}(T,0)^{k+1}x_n+
\lma_n  R^{(\mu_n,\lma_n)}(T,0)^{k} \left(\int_{0}^{T} h_n(s) \d s
\right),
\end{eqnarray*}
where $h_n(s):=R^{(\mu_n,\lma_n)}(T,s) F(s,\Psi_{s}^{(\lma_n)}(x_n,\mu_n),\mu_n)$ for $s\in [0,T]$ and
$n\geq 1$. Let $\eps > 0$ be arbitrary and a sequence $(k_n)$ of positive integers be such that
$k_n \to +\infty$ and $k_n \lma_n\to \eps$ as $n\to +\infty$.  Then, reasoning as before, one obtains
\begin{equation}\label{rr2}
x_n = R^{(\mu_n,\lma_n)}(T,0)^{k_n}x_n +J_n \left(\int_{0}^{T}
h_n(s) \d s\right) \quad \mbox{ for \ }n\ge 1,
\end{equation}
where $J_n:=\lma_n\sum_{k=0}^{k_n-1} R^{(\mu_n,\lma_n)}(T,0)^{k}$.
Note that, in view of Lemma \ref{14092008-1339},
\begin{align*}
h_n(s)\to F(s,x_0,\mu_0) & \quad \mbox{ as } n\to+\infty, \mbox{ uniformly for }  s\in [0,T], \\
R^{(\mu_n,\lma_n)}(T,0)^{k_n}x_n \to S_{\widehat{A^{(\mu_0)}}} (\eps T) x_0 & \quad\mbox{ as } n\to +\infty \mbox{ and }\\
J_n x\to \frac{1}{T}\int_{0}^{\eps T} S_{\widehat{A^{(\mu_0)}}}(s)
x\d s & \quad  \mbox{ as } n\to +\infty, \mbox{ for any }  x\in E.
\end{align*}
Thus, after passing in (\ref{rr2}) to the limit with $n\to +\infty$, one has
$$
x_0=S_{\widehat{A^{(\mu_0)}}} (\eps T) x_0  + \left[\frac{1}{T}\int_{0}^{\eps T} S_{\widehat{A^{(\mu_0)}}}(\tau)  d\tau\right]\left( \int_{0}^{T} F(s,x_0,\mu_0) \d s\right),
$$
which rewritten as
$$
-\frac{1}{\eps T}\left( S_{\widehat{A^{(\mu_0)}}} (\eps T) x_0 -
x_0\right) = \frac{1}{\eps T} \int_{0}^{\eps T}
S_{\widehat{A^{(\mu_0)}}}(\tau)\widehat{F}(x_0,\mu_0) \d\tau.
$$
Letting $\eps\to 0^+$ yields $-\widehat{A^{(\mu_0)}}x_0= \widehat
F(x_0,\mu_0)$, a contradiction completing the proof. \hfill
$\square$

\begin{Lem}\label{25092008-1149} {\em (see \cite[Proposition 4.3]{Cwiszewski-1})}
Let $F:E\to E$ be a completely continuous locally Lipschitz with sublinear growth and let $\Xi_t:E\to E$ be the translation along trajectories operator by time $t > 0$ for the equation $$\dot u(t) = - u(t) +F(u(t)), \quad t\in[0,T].$$
Then, for each $t > 0$, the mapping $\Xi_t$ is a $k$-set contraction and if an open bounded  $U\subset E$ is such that
$0\not\in (I-F)(\partial U)$, then there exists $t_0>0$ such that for any $t\in (0,t_0]$, $\Xi_t(x)\neq x$ and
$$
\deg(I-F, U) = \deg(I-\Xi_t, U).
$$
\end{Lem}

\noindent\textbf{Proof of Theorem \ref{18092008-2220}.}
First we reduce the proof to the case where the nonlinear perturbation is compact. By Remark \ref{19072009-2137} (b) we infer that the operator $\widehat{A}$ and mapping $\widehat F$ satisfy assumptions of Lemma \ref{22092008-1201}. Therefore there is locally Lipschitz compact mapping $\widehat F_L :E \to E$ such that
\begin{equation}\label{13092009-2214}
\widehat A x + (1-\mu)\widehat F(x) + \mu \widehat F_L(x) \neq 0 \quad \mbox{ for \ } x\in
\partial U \cap D(\widehat A), \ \mu\in [0,1].
\end{equation}
Thus, applying Lemma \ref{19092008-0014} to equations associated
to
$$\dot u(t)=\lma A(t) u(t)+ \lma ((1 - \mu)\widehat F(u(t)) + \mu \widehat F_L(u(t))), \quad t\in [0,T], $$
and the homotopy invariance of the topological degree, provide \\
\noindent {\bf Claim A.} {\em There exists $\lma_1>0$ such that for any $\lma\in (0,\lma_1]$
$$
\deg(I-\Phi_{T}^{(\lma)}, U) = \deg(I-\widetilde\Phi_{T}^{(\lma)}, U),
$$
where $\widetilde \Phi_{T}^{(\lma)}:\overline U\to E$ is the
translation along trajectories operator by the time $T$ for the
equation
$$
\dot u(t) = \lma A(t)u(t) + \lma \widehat F_L (u(t)).
$$}
\indent Next we prove

\noindent {\bf Claim B.} {\em There exists $\lma_2\in (0, \lma_1]$ such that for any $\lma\in (0,\lma_2]$
\begin{equation} \label{25092008-1143}
\deg(I-\widetilde\Phi_{T}^{(\lma)}, U) = \deg(I-\overline\Phi_{T}^{(\lma)}, U),
\end{equation}
where $\overline \Phi_{T}^{(\lma)}$ is the translation along trajectories operator
by the time $T$ for the equation
$$\dot u(t) = - \lma u(t) - \lma \widehat A^{-1}\widehat F_L(u(t)), \quad t\in [0,T].$$}
To this end, consider a differential problem given by
\begin{equation}\label{22092008-2350}
\dot u(t) = \lma \widetilde A^{(\mu)}(t)u(t) + \lma \widetilde F(u(t),\mu) \quad\mbox{ on } [0,T],
\end{equation}
where
\begin{align*}
\widetilde A^{(\mu)}(t) & := -\mu I + (1-\mu) A(t) &&\hspace{-20mm} \mbox{ for \ } t\in [0,T], \\
\widetilde F(x,\mu) & := [(1-\mu)I - \mu \widehat A^{-1}]\widehat F_L(x) &&\hspace{-20mm} \mbox{ for \ } x\in E, \ \mu\in [0,1].
\end{align*}
Lemma \ref{14092008-1339} shows that the family
$\{\lambda\widetilde A^{(\mu)}(t)\}_{t\in [0,T]}$ satisfies
$(Hyp'_1)$, $(Hyp_2)$ -- $(Hyp_4)$ and so the family $\{\widetilde
A^{(\mu)}(t)\}_{t\in [0,T]}$ fulfills the assumptions of Lemma
\ref{19092008-0014}. It is also clear that $\widetilde F$ is
locally Lipschitz in $x$ uniformly with respect to $\mu$ and
compact, which, in particular, means that it has sublinear growth
uniformly with respect to $\mu$. For any $\lambda\in(0,\infty)$,
let $\Psi:\overline U\times [0,1]\to E$ be given by
$$
\Psi_{T}^{(\lma)}(x,\mu):=u(T;x,\mu,\lma), \quad\mbox{ for \ }x\in E, \ \mu\in[0,1],
$$
where $u(\, \cdot\, ;x,\mu,\lma)$ stands for the mild solution of (\ref{22092008-2350}) starting at $x$.\\
\indent Observe that
\begin{equation}\label{22092008-2359}
[(1 - \mu) \widehat A - \mu I]x + \widetilde F(x,\mu)\neq 0 \quad\mbox{ for \ }
\mu\in [0,1],\, x\in \partial U\cap D(\widetilde A^{(\mu)}).
\end{equation}
Indeed, suppose that for some $\mu\in [0,1]$ and $x\in \partial U\cap D(\widehat A)$ we have
$[(1-\widetilde \mu) \widehat A -\mu I]x + \widetilde F(x,\mu)= 0.$
If $\mu=1$ then $-x-\widehat A^{-1}\widehat F_L(x)=0$, which contradicts (\ref{13092009-2214})
and if $\mu\in [0,1)$ then
$$
x=-R(\widehat A;\mu/(1-\mu)) (I-\mu/(1-\mu)\widehat A^{-1})\widehat F_L(x),
$$
which due to the resolvent identity gives $x=-\widehat{A}^{-1} \widehat F_L(x)$,
again a contradiction proving (\ref{22092008-2359}). Thus, applying Lemma \ref{19092008-0014}
and the homotopy invariance of the topological degree to $\Psi_{T}^{(\lma)}$, we find $\lma_2\in (0, \lma_1]$ such that
for any $\lma\in (0,\lma_2]$, (\ref{25092008-1143}) holds, which ends the proof of Claim B.\\
\indent Finally, by applying Lemma \ref{25092008-1149}, one gets $\lma_0\in (0,\lma_2]$ such that for any $\lma\in (0,\lma_0]$
\begin{equation*}
\deg(I+ \widehat {A}^{-1} \widehat F_L, U) = \deg(I-\overline \Phi_{\lma T}^{(1)},U).
\end{equation*}
Combining this with (\ref{13092009-2214}), we infer that, for
$\lambda\in(0,\lambda_0]$,
\begin{align*}
\Deg(\widehat A+ \widehat F,U) & = \Deg(\widehat A + \widehat F_L, U)
= \deg(I+ \widehat {A}^{-1} \widehat F_L, U)\\
& = \deg(I-\overline \Phi_{\lma T}^{(1)},U) = \deg(I-\overline \Phi_{T}^{(\lma)}, U),
\end{align*}
which together with Claims A and B completes the proof.\hfill $\square$

As an immediate consequence of Theorem \ref{18092008-2220}, we get the following result.
\begin{Cor}
If $0\not\in (\widehat A+\widehat F) (\partial U \cap D(\widehat A))$
and $\deg(\widehat A+\widehat F, U)\neq 0$, then $(P_{T,\lambda})$
admits a solution for small $\lma>0$.
\end{Cor}
By means of {\em a priori bounds} type assumption, we get the existence criterion for periodic solutions.
\begin{Th} {\em (Continuation principle)} \label{30092008-1514}
Let  a family $\{A(t)\}_{t\in [0,T]}$ and a mapping $F:[0,T]\times E\to E$ satisfy $(Hyp'_1)$, $(Hyp_2)$ -- $(Hyp_5)$ and $(F_1)$ -- $(F_3)$, respectively. If $(P_{T,\lambda})$ has no $T$-periodic points in $\partial U\times (0,1)$ and $\Deg(\widehat{A} + \widehat{F}, U)\neq 0$, then $(P)$ admits a mild solution $u:[0,T]\to E$ such that $u(0) = u(T)\in\overline U$.
\end{Th}
\noindent\textbf{Proof.} If $\Phi_T^{(1)}(x) = x$ for some $x\in\partial
U$, then the assertion holds. Hence, assume that $\Phi_T^{(1)}(x)\neq x$
for $x\in\partial U$. By Theorem \ref{18092008-2220}, there
exists $\lambda_0 \in (0,1)$ such that, for any $\lma\in (0,\lma_0]$, $\Phi_T^{(\lambda)}(x)\neq x$ and
\begin{equation}\label{22092008-2340}
\deg (I - \Phi_T^{(\lambda)}, U) = \Deg ( \widehat{A} + \widehat{F}, U).
\end{equation}
Then the mapping $\overline U\times [\lambda_0,1]\ni (x,\lambda)\mapsto \Phi_T^{(\lambda)}(x)$
provides an admissible homotopy (in the degree theory of $k$-set contraction vector fields)
and by the homotopy invariance
$$
\deg(I - \Phi_T^{(1)}, U) = \deg(I - \Phi_T^{(\lambda_0)},U),
$$
which together with (\ref{22092008-2340}) and the assumption implies the existence of
$x\in U$ such that $\Phi_{T}^{(1)}(x)=x$.\hfill $\square$

The above Continuation Principle can be useful when studying asymptotically linear evolution systems.
\begin{Th}\label{24102008-0927}
Let a family $\{A(t)\}_{t\in [0,T]}$ satisfy $(Hyp'_1)$ and $(Hyp_2)$ -- $(Hyp_5)$
and let $F:[0,T]\times E\to E$ be a completely continuous mapping satisfying $(F_1)$, $(F_2)$ and $(F_4)$. Assume also that $\{F_\infty (t):E\to E \}_{t\in [0,T]}$ is a family of compact linear operators such that the mapping $t\mapsto F_\infty (t)\in {\cal L} (E,E)$ is continuous on $[0,T]$, $F_\infty (0) = F_\infty (T)$ and
\begin{equation}\label{13072009-2241}
\lim_{\|x\|\to +\infty} \frac{\| F(t,x)-F_\infty(t)x\|}{\|x\|} = 0 \quad\mbox{ uniformly with respect to } t\in [0,T].
\end{equation}
If, for each $\lambda \in (0,1]$, the parameterized linear periodic problem
\begin{equation}\label{23012009-1323}
\left\{ \begin{array}{l}
\dot u(t)=\lambda (A(t)+F_\infty (t))u(t), \quad t\in [0,T]\\
u(0)=u(T)
\end{array}\right.
\end{equation}
has no nontrivial solution and $\Ker ({\widehat A}+ {\widehat F_\infty }) =\{ 0 \}$, then $(P)$ admits a $T$-periodic mild solution.
\end{Th}
\textbf{Proof.} We begin with proving that there exists $R_1>0$ such that
$0\not\in (\widehat A + \widehat F) ((E\setminus  B(0,R_1)) \cap D(\widehat A))$ and
\begin{equation}\label{26102008-2220}
\left|\Deg(\widehat A+\widehat F, B(0,R_1) ) \right| = 1.
\end{equation}
Define $H:[0,T] \times E \times [0,1] \to E$ by
$$
H(t,x,\lma):= \left\{\begin{array}{ll}
\lambda F(t,\lambda^{-1} x) & \mbox{ for \ } t\in [0,T], \ x\in E, \ \lambda\in (0,1],\\
                      F_\infty (t)x &  \mbox{ for \ } t\in [0,T], \ x\in E, \ \lambda=0.
\end{array}\right.
$$
Standard arguments show that both $H$ and $\widehat H$ are
completely continuous and for any $\lambda_0\in [0,1]$
\begin{equation}\label{26102008-2214}
\lim_{\|x\|\to + \infty,\, \lma\to \lma_0} \frac{\|\widehat H(x,\lma)-\widehat F_\infty x \|}{\|x\|}=0.
\end{equation}
Observe that there is $R_1>0$ such that
\begin{equation}
\widehat A x + \widehat H(x,\lambda)\neq 0
\ \ \mbox{ for \ } x\in (E\setminus B(0,R_1))\cap D(\widehat A), \ \lambda \in [0,1].
\end{equation}
Otherwise there are $(x_n)$ in $E$ and  $(\lma_n)$ in $[0,1]$ such that
$\widehat A x_n + \widehat H (x_n, \lma_n)=0$ and $\|x_n\|\to +\infty$.
Put $z_n:=x_n/\|x_n\|$. If $\lma_n=0$ for some $n\geq 1$, then
$z_n = - \widehat A^{-1} \widehat F_\infty z_n$,
a contradiction to the assumption. If $(\lma_n)$ in $(0,1]$, then
\begin{equation}\label{13072009-2237}
z_n= - \|x_n\|^{-1}{\widehat A}^{-1}\widehat H(x_n, \lma_n) = - {\widehat A}^{-1} (\lma_n^{-1}
\|x_n\|)^{-1}  \widehat F(\lma_{n}^{-1} \|x_n\| z_n).
\end{equation}
Since $\lim_{\|x\|\to +\infty} \|\widehat F(x)-\widehat F_\infty x\|/\|x\| = 0$ and $\rho_n := \lma_n^{-1} \|x_n\| \to +\infty$ as $n\to\infty$,
\begin{equation}\label{20072009-0116}
\|z_n + \widehat A^{-1} \widehat F_\infty z_n\| \le \|\widehat A^{-1}\|
\|\widehat F(\rho_n z_n) -  \widehat F_\infty (\rho_n z_n)\| / \rho_n \to \infty
\quad \mbox{ as \ }n\to \infty.
\end{equation}
By Lemma \ref{19092008-0029} the linear operator $\widehat
F_\infty$ is compact, which together with (\ref{20072009-0116})
means that $(z_n)$ contains a convergent subsequence. Hence, we
may assume that $z_n\to z_0$ for some $z_0\in E$ and by
(\ref{20072009-0116}) we infer that $z_0= -\widehat A^{-1}
\widehat F_\infty z_0$, which is again a contradiction
meaning that (\ref{26102008-2220}) holds for sufficiently large $R_1>0$.\\
\indent Now we claim that
\begin{equation} \label{23012009-1328}\mbox{\parbox[t]{130mm}{
there is $R\geq R_1$ such that, for any $\lma\in (0,1)$, the problem $(P_{T,\lma})$
has no periodic  solutions starting at points from $\partial B(0,R)$.} }
\end{equation}
Otherwise there exist
$(u_n)$ and $(\lma_n)$ in  $(0,1)$ such that, for each $n\geq 1$, $u_n$ is a solution of $(P_{T, \lma_n})$ and $\|u_n(0)\|  \to +\infty$ as $n\to +\infty$. Putting $v_n:=u_n / \|u_n\|_{\infty}$, where
$\|u_n\|_\infty:=\max_{t\in [0,T]} \|u_n(t)\|$, one has
\begin{equation}\label{25012009-2112}
v_n (t)= R^{(\lma_n)}(t,0) v_n(0) + \lma_n \|u_n\|_{\infty}^{-1}
\int_{0}^{t} R^{(\lma_n)} (t,s) F ( s , \|u_n\|_{\infty} v_n (s) ) \d s.
\end{equation}
Note that, by (\ref{13072009-2241}) and the fact that $F$ is a completely continuous mapping,
for any $\eps>0$, there is $m_\eps\geq 0$ such that
$\| F(t,x)-F_\infty(t)x\| \leq \eps\|x\| + m_\eps$ for  $x\in E$.
Consequently, for each $\eps > 0$, there exists $n_\eps\geq 1$ such that, for any $n\geq n_\eps$ and $s\in [0,T]$,
\begin{equation}\label{05022009-1319}
\ \ \left\| \|u_n\|_\infty^{-1} F(s, \|u_n\|_{\infty} v_n(s))- F_\infty(s) v_n(s) \right\| \leq \|u_n\|_{\infty}^{-1} (\eps \| u_n \|_\infty \|v_n(s)\| + m_\eps) \leq 2\eps.
\end{equation}
For each $n\geq 1$, put $h_n(s):= \|u_n\|_{\infty}^{-1} F ( s , \|u_n\|_{\infty} v_n (s) )$ for $s\in[0,T]$. If $s\in [0,T]$, then using (\ref{05022009-1319}) and the compactness of $F_\infty(s)$, for arbitrary $\eps > 0$, we deduce that
$$
\beta(\{ h_n(s) \}_{n\ge 1})\le \beta(F_\infty(s)(\{v_n(s)\}_{n\ge 1}))  + 2\eps = 2\eps\quad\mbox{ for \ }s\in[0,T],
$$
which, by passing to the limit with $\varepsilon\to 0$, gives $\beta(\{ h_n(s) \}_{n\ge 1}) = 0$. Obviously, we may also assume that $\lma_n\to\lma_0$ for some $\lma_0\in [0,1]$. \\
\indent If $\lma_0=0$, then note that
\begin{equation}\label{20072009-0255}
v_n(0)=R^{(\lma_n)}(T,0)^{k_n} v_n(0) + \left[\lma_n \sum_{k=0}^{k_n-1} R^{(\lma_n)} (T,0)^{k}\right] \left( \int_{0}^{T} R^{(\lma_n)}(T,s)h_n(s) \d s\right),
\end{equation}
where $(k_n)$ is an arbitrary sequence of positive integers.
If we put $k_n:=[T/\lma_n]$ for $n\geq 1$, then
$k_n\lambda_n \to T$ as $n\to +\infty$ and, in view of  Theorem \ref{29082008-1821} and Lemma \ref{19092008-0012}, we get
\begin{eqnarray*}
\beta(\{ v_n(0) \}_{n\ge 1}) \leq e^{-\omega T}\beta(\{ v_n(0) \}_{n\ge 1})
+ \frac{1-e^{-\omega T}}{\omega} \beta\left(\left\{ \int_{0}^{T} R^{(\lma_n)}(T,s)h_n(s) \d s\right\}_{n\ge 1}\right) \\
\le e^{-\omega T}\beta(\{ v_n(0) \}_{n\ge 1}) +
(1-e^{-\omega T})\omega^{-1}\int_{0}^{T} \beta(\{ h_n(s) \}_{n\ge 1}) \d s
= e^{-\omega T}\beta(\{ v_n(0) \}_{n\ge 1}).
\end{eqnarray*}
In consequence $\beta(\{v_n(0)\}_{n\ge 1})=0$. Furthermore, by (\ref{25012009-2112}) and Lemma \ref{19092008-0012}, for any $t\in [0,T]$, one has
\begin{align*}
\beta(\{v_n (t)\}_{n\ge 1}) & \le  \beta(\{R^{(\lma_n)}(t,0) v_n(0)\}_{n\ge 1})
+ \lma_n \int_{0}^{t} \beta(\{R^{(\lma_n)} (t,s) h_n(s)\}_{n\ge 1}) \d s \\
& \le \beta(\{v_n(0)\}_{n\ge 1}) +
\int_{0}^{t} \beta\left(\left\{\lma_n h_n(s)\right\}_{n\ge 1}\right) ds = 0,
\end{align*}
i.e. $\beta(\{ v_n (t)\}_{n\ge 1})=0$ for $t\in[0,T]$.
Since, by (\ref{05022009-1319}), the set $\{h_n\}_{n\ge 1}$ is bounded in $C([0,T],E)$, applying Proposition \ref{30122008-1629} (ii), we infer that $\{ v_n\}_{n\ge 1}$ is relatively compact in
$C([0,T],E)$ and without loss of generality, we assume that $v_n\to v_0$ in $C([0,T],E)$
and $v_n(0) \to x_0:=v_0 (0)$.
Furthermore, by (\ref{05022009-1319}), for arbitrary $\eps > 0$, there exists $n_\eps \geq 1$ such that, for any $n\geq n_\eps$ and $t\in [0,T]$,
\begin{align*}\label{20072009-0249}
\|v_n (t) - R^{(\lma_n)}(t,0) v_n(0)\| & \leq
\lma_n \int_{0}^{t} \|R^{(\lma_n)} (t,s) F_\infty(s) v_n (s)\| \d s+ 2\lma_n \eps T \\
& \leq  \lma_n K \int_{0}^{t} \|v_n (s)\| \d s +2\lma_n \eps T \to 0 \quad\mbox{ as }n\to \infty,
\end{align*}
where $K:=\sup_{\tau\in[0,T]}\|F_\infty(\tau)\|$. This together with Lemma \ref{14092008-1339} (i) imply that $v_0(t) = x_0$ for any $t\in E$ and in particular $x_0 \neq 0$, since $\|v_0\|_\infty = 1$. Hence, in view of (\ref{05022009-1319}) we find that
\begin{equation}\label{20072009-0305}
h_n(s) \to F_{\infty} (s)x_0 \quad \mbox{ as }n\to \infty \mbox{ uniformly for } s\in [0,T].
\end{equation}
Now fix an arbitrary $\eps>0$ and take any sequence $(k_n)$ of positive integers such that
$k_n \lma_n \to \eps$.
Applying Lemma \ref{14092008-1339} (ii), (iii) and (\ref{20072009-0305}) and passing to the limits in (\ref{20072009-0255}), we obtain
$$
x_0 = S_{\widehat A} (\eps T) x_0 + \left[ \frac{1}{\eps T} \int_{0}^{\eps T} S_{\widehat A}(\tau) d\tau \right] \left( \int_{0}^{T} F_\infty(s) x_0 \d s \right),
$$ i.e.
$$
-\frac{1}{\eps T} \left(S_{\widehat A} (\eps A) x_0 - x_0\right) = \frac{1}{\eps T}\int_{0}^{\eps T} S_{\widehat A }(\tau) \widehat F_\infty x_0 d\tau,
$$
Hence, a passing to the limit with $\eps\to 0$ yields $x_0\in D(\widehat A)$ and $(\widehat A + \widehat F_\infty)x_0=0$, a contradiction proving (\ref{23012009-1328}) in the case $\lma_0 = 0$.\\
\indent If $\lma_0\in (0,1]$, then, by (\ref{25012009-2112}) and Lemma \ref{19092008-0012},
\begin{align*}
& \beta(\{v_n(0)\}_{n\ge 1}) \le  e^{-\lma_0 \omega T}\beta(\{v_n(0)\}_{n\ge 1}) +
\lambda_0\beta\left(\left\{ \int_{0}^{T} R^{(\lma_n)}(T,s)h_n(s) \d s\right\}_{n\ge 1}\right) \\
& \le e^{-\lma_0 \omega T}\beta(\{v_n(0)\}_{n\ge 1}) +
\int_{0}^{T} \lma_0 e^{-\lma_0 \omega (T-s)}\beta\left(\left\{ h_n(s)\right\}_{n\ge 1}\right) ds = e^{-\lma_0 \omega T}\beta(\{v_n(0)\}_{n\ge 1})
\end{align*}
and, consequently, $\beta(\{v_n(0)\}_{n\geq 1})=0$. Using again
(\ref{25012009-2112}), for all $t\in [0,T]$,
\begin{equation*}
\beta(\{v_n(t)\}_{n\ge 1}) \le  e^{-\lma_0 \omega t}\beta(\{v_n(0)\}_{n\ge 1}) +
\int_{0}^{t} \lma_0 e^{-\lma_0 \omega (t-s)}\beta\left(\left\{ h_n(s)\right\}_{n\ge 1}\right) ds = 0,
\end{equation*}
which gives $\beta(\{ v_n(t)\}_{n\ge 1})=0$ for any $t\in [0,T]$.
Hence, due to the boundedness of $\{h_n\}_{n\ge 1}$ in $C([0,T],E)$ and Proposition \ref{30122008-1629} (ii), it follows that $\{ v_n\}_{n\ge 1}$ is relatively compact in $C([0,T],E)$ and, without loss of generality, we assume that $v_n\to v_0$ in $C([0,T],E)$.
Then, using (\ref{25012009-2112}), (\ref{05022009-1319}) and Proposition \ref{14012008-1812},
we infer, that for any $t\in [0,T]$,
$$
v_0(t)= R^{(\lma_0)} (t,0)v_0 (0) + \lma_0 \int_{0}^{t} R^{(\lma_0)} (t,s) F_\infty (s) v_0(s) \d s,
$$
i.e. $v_0$ is a nontrivial mild solution of (\ref{23012009-1323}) with $\lma=\lma_0\in(0,1]$,
which is a contradiction proving (\ref{23012009-1328}). \\
\indent Finally, (\ref{26102008-2220}) and (\ref{23012009-1328})
allow us to apply Theorem \ref{30092008-1514} to finish the proof.
\hfill $\square$

\section{An application to hyperbolic partial differential equations}
We end the paper with an example of a periodic problem for the
hyperbolic evolution equation with a time-dependent damping term.
Suppose that $\Omega$ is an open bounded subset of $\R^N$ and
$A:D(A)\to E$ is a positive self-adjoint linear operator with
compact resolvents defined on a Hilbert space $X:=L^2(\Omega)$
with the scalar product and the corresponding norm denoted
by $(\cdot,\cdot)_0$ and $|\cdot|_0$, respectively. It is well known that
such $A$ determines its fractional power space $X^{1/2}$ being a Hilbert
space as well. If we denote the scalar product and the corresponding norm by $(\cdot,\cdot)_{1/2}$ and
$|\cdot|_{1/2}$, respectively, then it is known that
$$
|u|_{1/2} \ge \lambda_1^{1/2}|u|_0 \quad \mbox{ for any } u\in X^{1/2}
$$
where $\lma_1>0$ is the smallest eigenvalue of $A$. Typical examples of $A$ satisfying
these conditions is $-\Delta_D$, where $\Delta_D $ is the Laplacian operator with zero the Dirichlet
boundary conditions  or  $-\Delta_N + \alpha I$, where $\Delta_N$ is the Laplacian
operator with the zero Neumann boundary conditions and $\alpha >0$.\\
\indent Consider a periodic problem
\begin{equation}\label{08082009-0110}
\left\{
\begin{array}{ll}
u_{tt} (x,t) + \beta (t)u_t (x,t) + (A u)(x, t) + f(t,u(x,t)) = 0  & \quad \mbox{ in } \Omega\times(0,T] \\
u(x,0)=u(x,T), \ u_t (x,0) = u_t (x,T) & \quad \mbox{ on }
\partial\Omega,
\end{array}
\right.
\end{equation}
where $\beta:[0,T]\to \R$ is a $T$-periodic continuously
differentiable function such that $\beta(t)>0$, for $t\in [0,T]$ and
$f:[0,T]\times \R \to \R$ is a continuous function satisfying the following properties
\begin{align}\label{z1}
& \mbox{there is } L>0 \mbox{ such that }  |f(t,s_1)\! -\! f(t,s_2)| \leq L|s_1 \!-\! s_2| \mbox{ for \ }t\!\in\! [0,T], \ s_1, s_2\in\mathbb{R}, \\ \label{z2}
& \mbox{there is } c>0 \mbox{ such that } |f(t,s)| \leq c (1+|s|) \mbox{ for \ } t\in[0,T], \ s\in\mathbb{R},\\ \label{z3}
& f(0,s)=f(T,s) \mbox{ for \ } s\in\mathbb{R}, \\ \label{z4}
& \lim_{|s| \to +\infty} \displaystyle{\frac{f(t, s)}{s}} =
f_{\infty} \mbox{ uniformly with respect to } t\in [0,T],
\end{align}
for some $f_{\infty}\in \R\setminus \sigma (A)$. If we define
$N_f:[0,T] \times X\to X$ by $N_f (t,u)(x):= f(t,u(x))$ for a.e.
$x\in \Omega$ and $t\in[0,T]$, then (\ref{08082009-0110}) can be
rewritten as a system
$$
\left\{
\begin{array}{ll}
\dot u(t)  =   v(t)\\
\dot v(t)  = - A u(t) - \beta (t) v(t) - N_f (t,u(t)), \quad\mbox{ for \ }t\in [0,T]
\end{array}
\right.
$$
and in a matrix form as
\begin{equation}\label{30062009-1633}
\dot z(t) = {\bold A} (t) z(t)+ {\bold F} (t,z(t)), \quad\mbox{ for \ }t\in [0,T]
\end{equation}
with operators ${\bold A}(t):D({\bold A}(t))\to {\bold E}$, $t\in [0,T]$, on the separable Banach space ${\bold E}:=X^{1/2} \times X$,
defined by
\begin{align}\label{06082009-1922}
D({\bold A} (t)) & := X^{1}\times X^{1/2} && \hspace{-15mm}\mbox{ for \ } t\in [0,T],\\
{\bold A} (t) (u, v) & :=(v,-Au-\beta (t)v) && \hspace{-15mm}\mbox{ for \ } t\in [0,T], \ (u, v) \in D({\bold A} (t)) \label{06082009-1923}
\end{align}
and ${\bold F}:[0,T] \times {\bold E}\to {\bold E}$ given by
${\bold F}\left(t, (u, v )\right) := (0, -N_f (t,u) )$ for $t\in [0,T]$, $(u, v) \in \textbf{E}$. \\
\indent  We claim that the family $\{ {\bold A}(t)\}_{t\in [0,T]}$ and the map ${\bold F}$ satisfy the assumptions of
Theorem \ref{24102008-0927} provided ${\bold E}$ is endowed with a proper norm.
To this end, for $\eta > 0$, define a new scalar product
$(\cdot,\cdot)_{{\bold E},\eta}:{\bold E}\times {\bold E} \to\R$, by
$$
\left( (u_1,v_1), (u_2, v_2)\right)_{{\bold E},\eta}: = (u_1,u_2)_{1/2} + (v_1+\eta u_1, v_2 +\eta u_2)_0.
$$
Clearly it is a well-defined scalar product and the corresponding
norm $\| \cdot \|_{{\bold E},\eta}$ is equivalent to the usual
product norm $\|\cdot\|$ in ${\bold E}=X^{1/2}\times X$. Let
$\beta_0 > 0$ be such that $\beta(t)\ge \beta_0$ for $t\in [0,T]$.
Putting $\gamma:=\max_{t\in [0,T]} \lambda_{1}^{-1/2}(\beta(t) +
1)$ for $0<\eta\leq 1$, one has
\begin{align*}
\left( {\bold A}(t)(u,v),(u, v)\right)_{{\bold E},\eta}
& = (v,u)_{1/2}+(- A u - \beta(t)v+\eta v, v + \eta u)_0 \\
& \hspace{-35mm} = (v,u)_{1/2} -(Au,v)_0 - (\beta(t)v,v)_0 + \eta|v|_{0}^2 - \eta (Au,u)_0 - \eta(\beta(t)v,u)_0 + \eta^2 (v,u)_0\\
& \hspace{-35mm} \leq -\eta|u|_{1/2}^{2} - (\beta_0-\eta)|v|_0^2 + \eta(\beta(t) + 1)|(v,u)|_0 \\
& \hspace{-35mm}\leq -\eta|u|_{1/2}^{2} - (\beta_0 - \eta)|v|_0^2 + \eta\gamma |u|_{1/2}|v| \\
& \hspace{-35mm}\leq  -\eta|u|_{1/2}^{2} - (\beta_0 - \eta)|v|_0^2 + (\eta/2) |u|_{1/2}^{2}+ (\eta\gamma^2/2)|v|_{0}^{2}\\
& \hspace{-35mm} = -(\eta/2)|u|_{1/2}^{2} -(\beta_0 - \eta - \eta \gamma^2/2)|v|_{0}^{2},
\end{align*}
and therefore, decreasing $\eta>0$ if necessary, there exists
$\omega=\omega(\eta)>0$ such that, $({\bold
A}(t)(u,v),(u,v))_{{\bold E},\eta} \leq -\omega \|(u,v)\|_{{\bold
E},\eta}^{2}$ for any $(u,v)\in D({\bold A} (t)) = X^{1}\times
X^{1/2}$. Since $\Im \, {\bold A}(t) = {\bold E}$, for $t\in
[0,T]$, we infer that, for $t\in [0,T]$, the operator ${\bold
A}(t)$ is a generator of $C_0$ semigroup satisfying
$$
\|S_{{\bold A}(t)}(s)\|_{{\bold E},\eta} \le e^{-\omega s}
\quad\mbox{ for } s \ge 0
$$
and, in particular, condition $(Hyp'_1)$ holds.
Moreover, observe that, for each $(u,v)\in X^{1}\times X^{1/2}$, the map $t \mapsto {\bold A}(t)(u,v)\in {\bold E}$ is continuously differentiable on $[0,T]$ as $\beta$ is so. Hence, in view of Proposition \ref{27022009-1214}, the family
$\{{\bold A}(t)\}_{t\geq 0}$ satisfies also conditions $(Hyp_2)$ and $(Hyp_3)$ with
${\bold V}:=X^{1}\times X^{1/2}$ equipped with
the norm given by $\|(u,v)\|_{{\bold V}}:= \|{\bold A(0)}(u,v)\|_{{\bold E},\eta}
+ \|(u,v)\|_{{\bold E},\eta}$ for $(u,v)\in {\bold V}$. By the periodicity of $\beta$, $(Hyp_5)$ holds.
Furthermore, observe that for $(u,v)\in {\bold V}=X^{1}\times X^{1/2}$
$$
{\bold A}_0 (u,v):= \frac{1}{T}\int_{0}^{T} {\bold A}(\tau)(u,v) d \tau = (v,-A u - \widehat \beta v)
$$
where $\widehat\beta:=(1/T)\int_{0}^{T}\beta(\tau)d\tau$ and $\Im \, {\bold A}_0 = {\bold E}$ which implies $(Hyp_4)$.
It can be easily verified that ${\bold A}_0$ is closed and consequently, $\widehat {\bold A} = {\bold A}_0$.
Since $\beta$ is $T$-periodic function, we conclude that $(Hyp_5)$ is also satisfied.\\
\indent It may be checked that ${\bold F}$ is continuous and, by
(\ref{z1}), (\ref{z2}) and (\ref{z2}), satisfies conditions
$(F_1)$, $(F_2)$ and $(F_4)$. Since the operator $A$ has compact
resolvents, the inclusion $X^{1/2}\subset X$ is compact and
therefore, both ${\bold F}$ and ${\bold F}_\infty:{\bold E}\to
{\bold E}$, given by ${\bold F}_{\infty} (u,v):=(0, - f_\infty
u)$, are completely continuous. Furthermore observe that
\begin{equation}\label{30072009-2032}
\limsup_{\|(u,v) \|\to \infty,\, t\to t_0} \frac{\|{\bold F}(t,(u,v)) -
{\bold F}_\infty (u,v)\|}{\|(u,v)\|} \leq \limsup_{|u|_{1/2}\to
\infty,\, t\to t_0} \frac{|N_f(t,u)-f_{\infty}u|_{0}}{|u|_{1/2}}.
\end{equation}
Now suppose that $(u_n)$ is a sequence in $X^{1/2}$ such that
$|u_n|_{1/2} \to +\infty$ as $n\to+\infty$ and $(t_n)$ in $[0,T]$
is such that $t_n\to t_0$ as $n\to +\infty$. If we put
$z_n:=u_n/|u_n|_{1/2}$ for $n\geq 1$, then by the compactness of
the inclusion $X^{1/2}\subset X$, there exists subsequence
$(z_{n_k})$ and $z_0\in X$ such that $z_{n_k} \to z_0$ in $X$ as
$k\to +\infty$. Without lost of generality we may assume that
$z_{n_k}(x) \to z_0(x)$, a.e. on $\Omega$, and there is $g\in X$
such that, for each $k\ge 1$, $|z_{n_k}(x)| \leq g(x)$ a.e. on
$\Omega$. Then, putting $\mu_n:= |u_n|_{1/2}$, by (\ref{z4}), we
find that
$$
|\mu_{n_k}^{-1} f(t_{n_k},x,\mu_{n_k} z_{n_k}(x)) - f_{\infty}z_{n_k} (x)|^{2} \to 0 \quad
\mbox{ a.e. on }\Omega, \ \mbox{ as }k\to \infty$$
and, for each $k\geq 1$,
\begin{align*}
|\mu_{n_k}^{-1} f(t_{n_k},x,\mu_{n_k} z_{n_k}(x)) -
f_{\infty}z_{n_k} (x)|^{2} \le g_0(x) \quad \mbox{ a.e. on } \Omega
\end{align*}
with $g_0 := (c(m + g) + f_\infty g)^2$ where $m:=\sup \{
\mu_{n}^{-1} \mid n\geq 1 \}$. Since $\Omega$ is a bonded set,
$g_0$ is integrable and by the Lebesgue dominated convergence
theorem
$$
|N_f(t_{n_k},u_{n_k})-f_\infty u_{n_k}|_0/|u_{n_k}|_{1/2} =
\int_{\Omega} |\mu_{n_k}^{-1} f(t_{n_k},x,\mu_{n_k} z_{n_k}(x)) -
f_{\infty}z_{n_k} (x)|^{2} dx \to 0
$$
as $k\to +\infty$, which together with (\ref{30072009-2032}), implies that
$$
\lim_{\|z \|\to \infty,\, t\to t_0} \frac{\|{\bold F}(t,z) - {\bold
F}_\infty z\|}{\|z\|} = 0
$$
and condition (\ref{13072009-2241}) is satisfied.\\
The following lemma will be helpful in verifying that (\ref{23012009-1323}) has no
nontrivial $T$-periodic solutions for $\lma\in (0,1)$.
\begin{Lem}\label{06082009-2014}
Let $A:D(A)\to E$ be a positive self-adjoint operator with compact
resolvents on a Hilbert space $X$ and, for fixed $\bar \beta>0$,
${\bold A}:D({\bold A})\to {\bold E}$ be an linear operator on ${\bold E}:= X^{1/2} \times X^{0}$ given by
$D({\bold A}):= X^{1}\times X^{1/2}$ and
$$
{\bold A}(u,v):= (v, - Au-\bar \beta v) \quad\mbox{ for \ } (u,v)\in D({\bold A}).
$$
Let ${\bold E}_k:= X_k \times X_k$, $k\geq 1$, where $X_k$ is the space spanned by the first $k$ eigenvectors of $A$ (corresponding to the first $k$ smallest eigenvalues of $A$), and ${\bold A}_k:{\bold E}_k\to {\bold E}$ be given by
$$
{\bold A}_k(u,v):={\bold A}(u,v) \quad\mbox{ for \ } (u,v)\in {\bold E}_k.
$$
Then\\
\pari{(i)}{${\bold A} ({\bold E_k})\subset {\bold E}_k$ for each $k\geq 1$;}\\
\pari{(ii)}{$R(\mu;{\bold A}_k)(u,v)= R(\mu;{\bold A}) (u,v)$ for any $k\geq 1$ and $(u,v)\in {\bold E}_k$ and $\mu>0$;}\\
\pari{(iii)}{$S_{{\bold A}_k}(t)(u,v) = S_{{\bold A}}(t)(u,v)$ for any $k\geq 1$ and $(u,v)\in {\bold E}_k$.}
\end{Lem}
{\bf Proof. } (i)  comes  straightforwardly from the fact that $A(X_k)\subset X_k$ for $k\geq 1$.\\
\indent (ii) If $(p,q)\in {\bold E}_k$ and $(u,v):=R(\mu;{\bold A}_k)(p,q)$ for some $\mu>0$, then
$\mu u - v = p$ and $Au+(\mu+\bar \beta)v=q$, i.e. $\mu(\mu+\bar \beta)u+Au = (\mu+\bar \beta)p+ q\in X_k$, which shows that
$u\in X_k$ and $v=\mu u-p\in X_k$. Therefore, $\mu(u,v)-{\bold A}_k(u,v)= \mu (u,v)-{\bold A}(u,v)=(p,q)$,
that is $R(\mu;{\bold A})(p,q)=(u,v)= R(\mu;{\bold A}_k)(p,q)$.\\
\indent (iii) follows from the Euler formula $S_{\bold A}(t)(u,v) = \lim_{n\to +\infty} (n/t)^{n}R(n/t; {\bold A})^{n} (u,v)$, $(u,v)\in {\bold E}_k$,
and (ii).  \hfill $\square$
\begin{Lem}\label{07082009-1207}
Let $\{ {\bold A}(t)\}_{t\in [0,T]}$ be given by {\em
(\ref{06082009-1922}), (\ref{06082009-1923})} and $\{{\bold A}_k
(t)\}_{t\in [0,T]}$, $k\geq 1$, be given by ${\bold A}_k (t)
(u,v):= {\bold A}(t)(u,v)$ for $(u,v)\in {\bold E}_k$. If
$\{{\bold R}(t,s) \}_{0\leq s \leq t\leq T}$ and $\{{\bold
R}^{(k)} (t,s)\}_{0\leq s\leq t\leq T}$, $k\geq 1$ are the
evolution systems determined by $\{ {\bold A}(t)\}_{t\in [0,T]}$
and $\{{\bold A}_k (t):{\bold E}_k\to {\bold E}_k\}_{t\in [0,T]}$,
$k\geq 1$, respectively, then, for any $k\geq 1$ and $(u,v)\in
{\bold E}_k$,
$$
{\bold R}^{(k)}(t,s) (u,v) = {\bold R}(t,s)(u,v) \quad\mbox{ for \ } 0\leq s\leq t\leq T.
$$
\end{Lem}
{\bf Proof}. By the construction of evolution systems (see \cite[Ch. 5, Theorem 3.1]{Pazy}), for any $(u,v)\in {\bold E}$,
\begin{equation}\label{06082009-2044}
{\bold R}(t,s) (u,v) = \lim_{n\to +\infty} {\bold R}_{n} (t,s) (u,v) \quad\mbox{ for \ } 0\leq s\leq t\leq T,
\end{equation}
where ${\bold R}_{n} (t,s):{\bold E}\to {\bold E}$, $n\geq 1$ are given by
$$ {\bold R}_{n} (t,s)\!:= \!\!\left\{
\begin{array}{ll}
\!\! S_{{\bold A}(t_j^{n})}(t\!-\!s) &  \mbox{ if } s,t\in [t_{j}^{n}, t_{j+1}^{n}], s \le t,\\
\!\! S_{{\bold A}(t_r^{n})}(t\!-\!t_{r}^{n})\!\left
(\prod\limits_{j=l+1}^{r-1} \!\!S_{{{\bold A}(t_j^{n})}}(T/n)\!
\right)\!S_{ {\bold A}(t_l^{n}) }( t_{l+1}^{n}\!-\!s)& \mbox{ if }
l < r \mbox{ and } s\in [t_{l}^{n}, t_{l+1}^{n}], \\ & \ \ \ \
t\in [t_{r}^{n}, t_{r+1}^{n}]
\end{array}
\right.$$
with $t_{j}^{n}:=(j/n)T$ for $j=0,1,\ldots, n$.
Similarly, for any $k\geq 1$ and $(u,v)\in {\bold E}_k$,
\begin{equation}\label{06082009-2045}
{\bold R}^{(k)} (t,s) (u,v) = \lim_{n\to +\infty} {\bold R}_{n}^{(k)} (t,s) (u,v) \quad\mbox{ for \ } 0\leq s\leq t\leq T.
\end{equation}
where
$$
{\bold R}_{n}^{(k)} (t,s)\!:= \!\!\left\{
\begin{array}{ll}
\!\! \!\!S_{{\bold A}_k (t_j^{n})}(t\!-\!s) &\!\!\!\!\mbox{if }  s,t\in [t_{j}^{n}, t_{j+1}^{n}], \, s\le t,\\
\!\!\! \!S_{{\bold A}_k
(t_r^{n})}(t\!-\!t_{r}^{n})\!\left(\prod\limits_{j=l+1}^{r-1}
\!\!S_{{{\bold A}_k(t_j^{n})}}(T/n)\! \right)\!S_{ {\bold
A}_k(t_l^{n}) }( t_{l+1}^{n}\!\!-\!s) & \!\!\!\!\mbox{if } l < r
\mbox{ and } s\in [t_{l}^{n}, t_{l+1}^{n}], \\ & \ t\in
[t_{r}^{n}, t_{r+1}^{n}].
\end{array}
\right.
$$
Lemma \ref{06082009-2014} (iii) states that $S_{{\bold
A(t_{n}^{j})}} (s) (u,v) = S_{{\bold A_k(t_{n}^{j})}} (s)(u,v)$
for $k\geq 1$, $(u,v)\in {\bold E}_k$, $s\geq 0$, $n\geq 1$ and
$j\in \{1,\ldots, n \}$. Hence, using the formulae
(\ref{06082009-2044}) and (\ref{06082009-2045}) completes the
proof. \hfill $\square$
\begin{Lem} {\em (cf. \cite{Lad}, \cite{Rybakowski})}\label{30062009-1545}
If $f:[0,T]\to X^{0}$ is continuous and
$(u,v):[0,T]\to X^{1/2}\times X^{0}$ is a mild solution of
$$
(u(t),v(t))' = {\bold A}(t) (u(t),v(t)) + (0,f(t)), \quad
t\in[0,T],
$$
then
\begin{eqnarray}\label{07082009-1309}
& & \frac{1}{2} \frac{d}{d t} |u(t)|^2_0   =  (u(t),v(t))_0 \quad  \mbox{ for \ } t\in [0,T], \\
\label{06082009-2221} & & \frac{1}{2} \frac{d}{d t} \left( |u
(t)|_{1/2}^{2} + |v(t)|_{0}^{2}\right)   = - \beta(t) |v
(t)|_{0}^{2} + (f(t),v(t))_0 \quad  \mbox{ for \ } t\in [0,T].
\end{eqnarray}
\end{Lem}
{\bf Proof}. Let $(u,v):[0,T]\to {\bold E}$ be a mild solution of
$$
\dot z (t) = {\bold A}(t) z(t) + (0, f(t)), \quad t\in [0,T].
$$
Put $(\overline u_k, \overline v_k):= (\widetilde P_k u(0),  P_k v(0))$ for $k\geq 1$,
where $\widetilde P_k:X^{1/2}\to X_k$ and $P_k:X^0 \to X_k$ are the orthogonal projections.
Furthermore, let $(\widetilde u_k, \widetilde v_k):[0,T]\to {\bold E}_k$ be the mild solution of
\begin{equation}\label{06082009-2212}
\left\{
\begin{array}{ll}
(\dot u(t),\dot v(t))={\bold A}_k(t)(u(t),v(t)) + (0, P_k f(t)), \quad t\in [0,T]\\
(u(0),v(0))= ( \overline u_k, \overline v_k )
\end{array}\right.
\end{equation}
and for each $k\geq 1$ define $( u_k, v_k):[0,T]\to {\bold E}$ by
$(u_k(t), v_k(t)):=(\widetilde u_k(t), \widetilde v_k(t))$ for
$t\in [0,T]$. Then
$$
(\widetilde u_k(t), \widetilde v_k(t))= {\bold R}^{(k)}(t,0)
(\overline u_k, \overline v_k) + \int_{0}^{t} {\bold R}^{(k)}
(t,s) (0,P_k f(s)) \d s \quad \mbox{ for } t\in [0,T],
$$
and by Lemma \ref{07082009-1207}, one gets
$$
(u_k(t), v_k(t)) = {\bold R} (t,0) (\overline u_k, \overline v_k)
+ \int_{0}^{t} {\bold R}(t,s) (0, P_k f(s))\d s \quad\mbox{ for }
t\in [0,T].
$$
Further, since $(\overline u_k, \overline v_k)\to (u(0),v(0))$ in
${\bold E}$ and $P_k f(t)\to f(t)$ in $X^0$ as $k\to +\infty$
uniformly for $t\in [0,T]$, by Proposition \ref{30122008-1629} (i)
we infer that $(u_k, v_k)\to (u,v)$ in $C([0,T], {\bold E})$.
Finally, treating (\ref{06082009-2212}) as a system of ordinary
differential equations with the family $\{{\bold
A}_k(t)\}_{t\in[0,T]}$ of bounded operators we see that
$(\widetilde u_k, \widetilde v_k)$ is, in particular, a classical
solution. Therefore, we obtain
\begin{eqnarray*}
& & (u_k(t),\dot u_k(t))_{1/2} =  (u_k(t), v_k(t))_{1/2}\\
& & (v_k(t), \dot v_k(t))_{0} =  - (u_k(t), v_k(t))_{1/2} -
\beta(t) |v_k(t)|_{0}^{2} + (v_k(t), P_k f(t))_0,
\end{eqnarray*}
for any $t\in [0,T]$ and, as a result,
$$
\frac{1}{2} \frac{d}{d t}(|u_k(t)|_{1/2}^{2} + |v_k(t)|_{0}^{2}) =
- \beta(t) |v_k(t)|_{0}^{2} + (v_k(t), P_k f(t))_0 \quad\mbox{ for
\ } t\in [0,T].
$$
Thus, we see that both, the functions $|u_k|_{1/2}^{2} + |v_k|_{0}^{2}$, $k\geq 1$ and their derivatives converge uniformly on $[0,T]$, which gives
(\ref{06082009-2221}). To see (\ref{07082009-1309}) observe that
\begin{equation}\label{07082009-1340}
\frac{1}{2} \frac{d}{d t}|u_k(t)|_{0}^{2} = (u_k(t),\dot u_k(t))_{0} =  (u_k(t), v_k(t))_{0} \quad\mbox{ for \ } t\in [0,T]
\end{equation}
and $u_k(t) \to u(t)$, $v_k(t) \to v(t)$ is a space $X^0$, as $k\to \infty$, uniformly with respect to $t\in[0,T]$. Hence we see that the functions $|u_k|_0^{2}$, $k\ge 1$, and their derivatives are convergent uniformly and, by (\ref{07082009-1340}), we are done.
\hfill $\square$

Now return to our considerations of (\ref{30062009-1633}) and  suppose that, for some, $\lma\in (0,1]$, $(u,v):[0,T]\to E$ is a $T$-periodic mild solution of
$$
(u,v)' = \lma({\bold A} (t) (u(t),v(t))+ {\bold F}_\infty (u(t),v(t))), \quad t\in[0,T].
$$
If view of Lemma \ref{30062009-1545} we get
$$
\frac{1}{2} \frac{d}{d t} \left(|u(t)|_{1/2}^{2}+|v (t)|_{0}^{2}
\right) = - \lma \beta(t) |v (t)|_{0}^{2} - \lma f_{\infty}
(u(t),v(t))_0
$$
and, after integrating and using (\ref{07082009-1309}), one has
$$
0 < \int_{0}^{T} \lma \beta(t)|v(t)|_{0}^{2} \d t = -
\frac{1}{2}\int_{0}^{T} \lambda f_{\infty} (|u(t)|_{0}^2)' \d t =
0,
$$
a contradiction proving that (\ref{23012009-1323}) has no nontrivial $T$-periodic solutions.\\
Finally, by a direct calculation, we see that $\Ker (\widehat {\bold A} + {\bold F}_\infty)=\{0\}$ since $f_\infty\not\in \sigma(A)$. Thus, in view of Theorem \ref{24102008-0927}, problem (\ref{08082009-0110}) admits a $T$-periodic solution in the sense that (\ref{30062009-1633}) has a $T$-periodic mild solution.

\end{document}